\documentclass{amsart}
\usepackage{wasysym}
\usepackage[T1]{fontenc}
\usepackage{mathtools, amssymb, graphicx, tcolorbox, extarrows} 
\usepackage{xcolor} 
\usepackage{hyperref}
\usepackage{float} 
\restylefloat{figure} 
\usepackage{parskip}

\usepackage{etoolbox}

\theoremstyle{definition}
\numberwithin{equation}{section}
\newtheorem{shared}{Dummy}[section]
\newtheorem{theorem}[shared]{Theorem}
\newtheorem{lemma}[shared]{Lemma}
\newtheorem{corollary}[shared]{Corollary}
\newtheorem{definition}[shared]{Definition}
\newtheorem{example}[shared]{Example}

\newtheorem{remark}[shared]{Remark}
\newtheorem*{example*}{Example}

\title{The Math Teaching Atlas: Trails, Anchors, and Compass in Action}
\thanks{A previous version of this paper by the author appeared in \emph{Journal of Humanistic Mathematics} 16(1) (2026), 344--364: \url{https://scholarship.claremont.edu/jhm/vol16/iss1/23/}.}
\author{Ivan Z. Feng}
\address{Department of Mathematics, University of Southern California, Los Angeles, CA 90089-2532, USA}
\email{ifeng@usc.edu}
\urladdr{https://dornsife.usc.edu/ivan/}
\date{March 31, 2026}
\begin{document}

\begin{abstract}
Mathematics is a mountain, but students need more than descriptions of the view: they need a trail they can actually walk. This paper presents the Math Teaching Atlas, a framework for mathematical exposition built around route units (single steps with explicit justifications), routeways (chains of route units), and roadmaps of multiple valid routeways to the same destination. New ideas are introduced through familiar anchors, while a mathematical compass (motivation) and driving simulations (concrete examples) help students extend known routeways and construct new ones. The paper also develops a route geometry on pedagogical route graphs, thereby beginning a mathematization of mathematics pedagogy itself. Taken together, the framework aims to make mathematical explanations more visible, navigable, and teachable.
\end{abstract}

\maketitle
\tableofcontents

\section{Clear Trails: Providing a Routeway, Not Descriptive Words}\label{Sec:trails}

Mathematics is logical, precise, and beautiful; yet modern math is often taught through verbose narratives that deviate from its logical core. Professors frequently explain concepts as if narrating a story rather than giving clear, step-by-step directions. Similarly, textbooks describe mathematics rather than explicitly showing its logical structure. This approach is problematic.

Math should be taught using clear, logical structure. Instead of extensive explanatory paragraphs, we should provide students with a \emph{routeway}: a transparent, step-by-step logical progression where connections are explicit. 

\begin{figure}[H]
    \centering
\includegraphics[width=\textwidth]{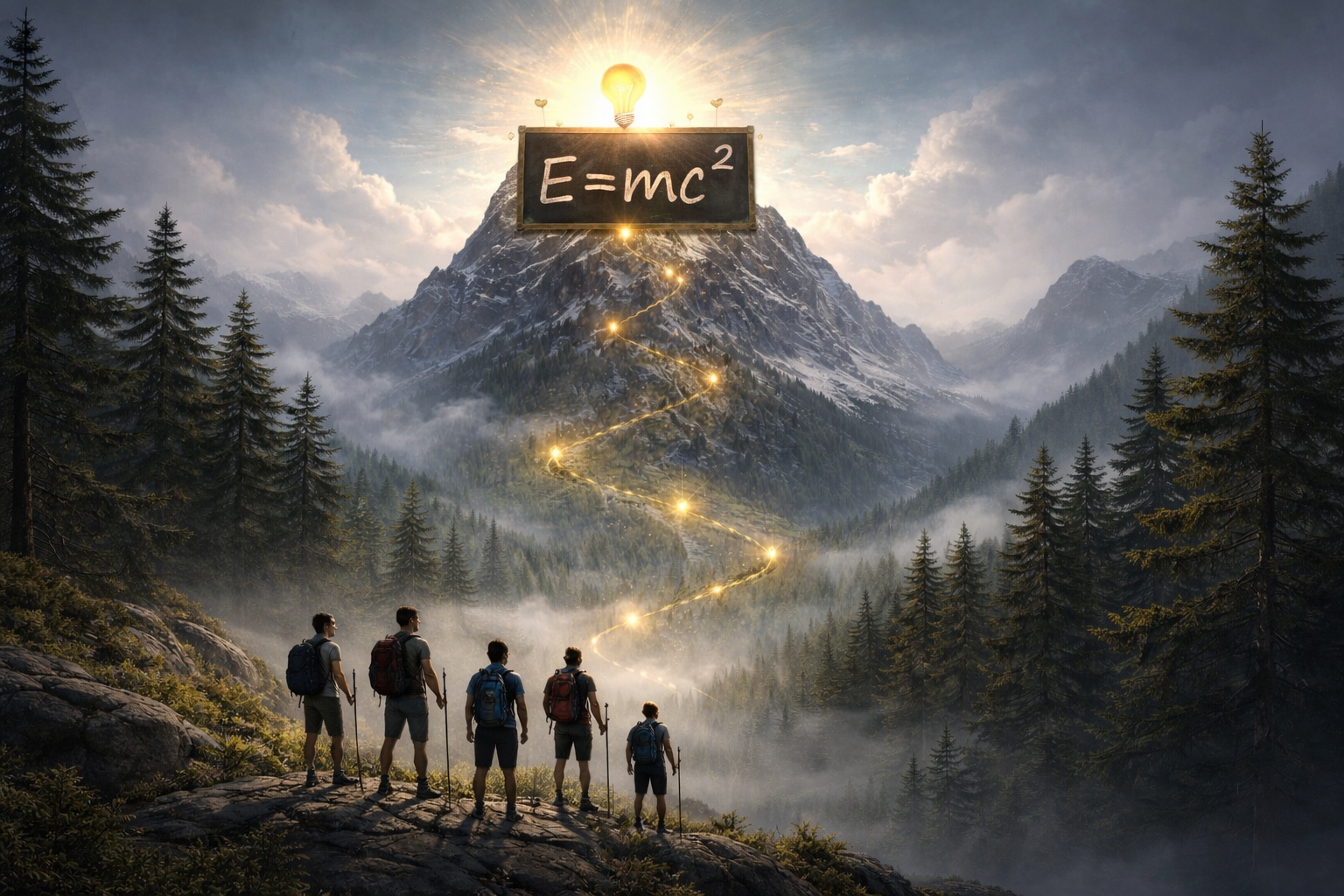}
    \caption{The summit alone is not enough: students need a visible trail.}
    \label{fig:mountain-trail}
\end{figure}

Consider a typical proof structure currently found in textbooks:

\begin{quote}
``Let’s prove Statement \( S \). Since \( A \), we have \( B \). From this, it follows that \( C \). Thus, \( D \) must hold. Therefore, we conclude \( E \).''\footnote{What's worse, \( E \) often does not match the original statement \( S \).}
\end{quote}

These statements merely announce facts rather than guide students through the logical connections. Students must mentally reconstruct missing steps and infer invisible links. 

Now contrast this with a clearly structured routeway, as illustrated in Figure \ref{fig:routeway}:

\begin{figure}[H]
    \centering
    \includegraphics[width=0.7\textwidth]{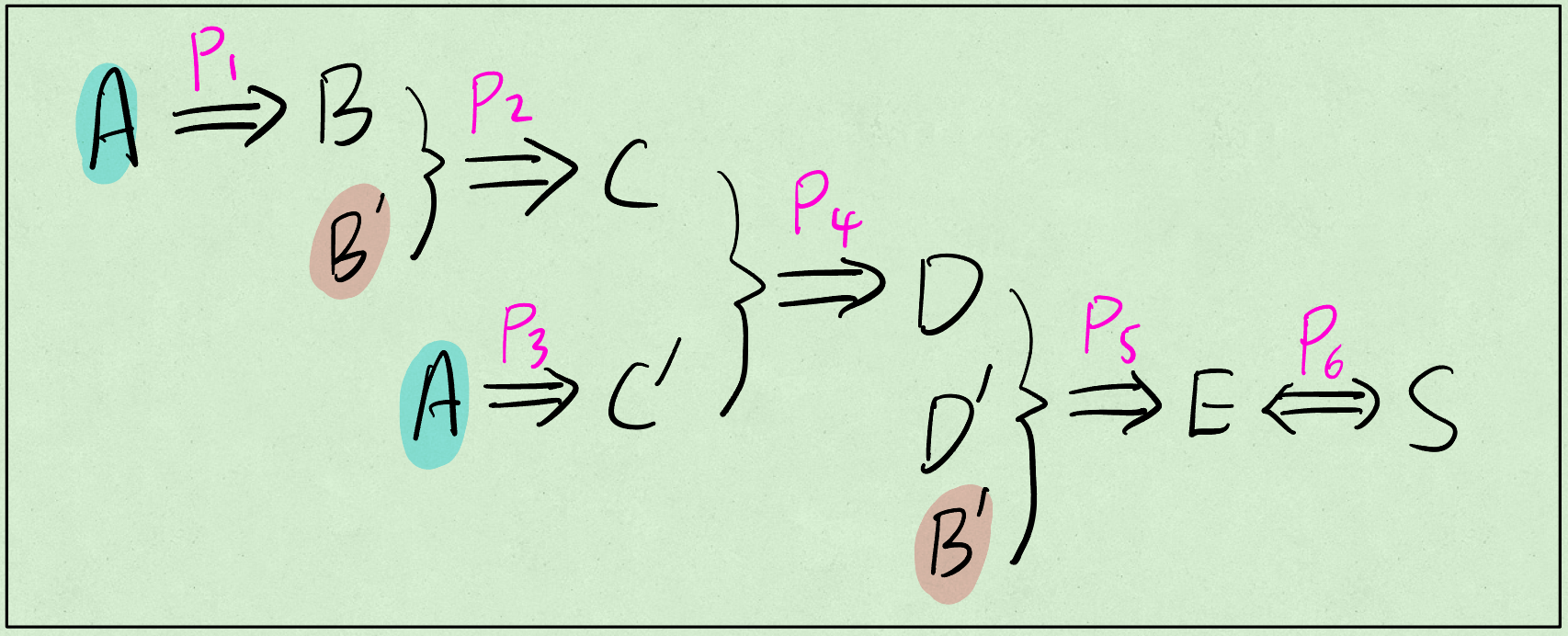}
    \caption{Routeway}
    \label{fig:routeway}
\end{figure}

As Lakatos argued in~\cite{lakatos}, proofs evolve through visible, refutable steps—not obscured by narrative leaps. In the structured approach above, every step is explicitly justified. If \( A \) leads to \( B \), it is because of the explicitly stated general fact \( P_1 \), which clearly states, ``Given \( A \), we have \( B \).'' If \( B \) then leads to \( C \), we understand precisely how this occurs: by combining another known condition \( B' \) (which is omitted in the textbook proof above) and then explicitly applying the justification \( P_2 \), which states ``\( B \) and \( B' \) imply \( C \).''

These general facts, denoted \( P_i \) (\( i=1,2,\dots \)), are often neglected in textbooks and lectures, buried in descriptive language. They may represent definitions, properties, or theorems, but their explicit statement bridges the logical gaps between steps. For instance, clearly stating \( P_1 \) as a general fact transforms the logical step from \( A \) to \( B \) into:
\[
A \ (\textbf{Premise}) \ \xRightarrow{P_1 \ (\textbf{Justification})} \ B \ (\textbf{Conclusion}).
\]

This premise--conclusion architecture itself goes back at least to Aristotle: in \textit{Prior Analytics} I.4, he writes, ``If \(A\) is predicated of all \(B\), and \(B\) of all \(C\), \(A\) must be predicated of all \(C\),'' \cite{aristotle-prior}. For example, a standard instance of this form is: every human is mortal; every Greek is human; therefore every Greek is mortal.

Even a student unfamiliar with mathematics can follow this reasoning: ``How did we get $B$ from $A$? Oh, I see—by $P_1$.'' Explicitly stating these justifications transforms gaps into bridges, uncertainty into confidence, and confusion into absolute clarity.

\subsection{Route Units}\label{Sec:route-units}

Let us formalize our discussion above into definitions using terminology consistent with the ``routeway'' theme:

\begin{definition}[Waypoint and route unit]
A \textbf{waypoint} is a single statement or a finite tuple of statements. In particular, a waypoint may have the form \(A=(A_1,\dots,A_m)\). For instance, in the route unit
\[
\begin{rcases*}
B\\
B'
\end{rcases*}
\xRightarrow{P_2} C,
\]
the waypoint is the pair \((B,B')\).

A \textbf{route unit} is a displayed inference of the form
\[
A \xRightarrow{P} B,
\]
where \(A\) and \(B\) are waypoints and \(P\), when present, is the \textbf{trail} explicitly given as the justification for the transition from \(A\) to \(B\). Here \(A\) is the \textbf{initial waypoint} and \(B\) is the \textbf{terminal waypoint}.

In particular, a \textbf{two-way route unit} is a route unit of the form
\[
A \xLeftrightarrow{P} B,
\]
where \(P\) explicitly justifies both directions \(A \Rightarrow B\) and \(B \Rightarrow A\).

A route unit is called \textbf{defective} if \(P\) is omitted.
\end{definition}

To formalize what counts as a ``single step,'' we fix the set of trails treated as atomic for the intended audience.

\begin{definition}[Pedagogical base field]\label{def:basefield}
For a fixed intended audience (course level, week, prerequisites, and conventions), the \textbf{(pedagogical) base field} \(\mathbb{B}\) is the set of trails (e.g., definitions, theorems, conventions, and allowed moves) treated as atomic, i.e., admissible as single-step justifications, for that audience.
\end{definition}

For example, in a textbook, we usually choose the pedagogical base field to consist of the material already covered in previous pages together with the relevant prerequisites.

With the pedagogical base field established, we can specify what constitutes a single, indivisible logical step.

\begin{definition}[$\mathbb{B}$-irreducible route unit]\label{def:B-irreducible}
Fix a pedagogical base field \(\mathbb{B}\). A non-defective route unit
\[
A \xRightarrow{P} B
\]
is \textbf{$\mathbb{B}$-irreducible} if \(P\in\mathbb{B}\) and \(B\) follows from \(A\) by a \emph{single application} of \(P\) (allowing only parameter-substitution/renaming, with no additional lemmas or intermediate route units). Equivalently, \(P\) produces \(B\) from the waypoint(s) \(A\) in one step: after possibly renaming symbols and specializing parameters, \(P\) reads ``from \(A\) infer \(B\).''
\end{definition}

Irreducibility depends on the chosen pedagogical base field \(\mathbb{B}\) in Definition~\ref{def:basefield}. 

\begin{example}\label{ex:ineq-basefield}
Consider the trail \(P\): ``Multiplying both sides of an inequality by a positive number \(c\) preserves the inequality.'' In an elementary base field, the route unit\footnote{Although commonly used interchangeably, this paper distinguishes carefully between ``\(=:\)'' and ``\(:=\)''. Specifically, the entity placed next to the colon (\(:\)) is always the one being newly introduced. For example, ``\(c =: 2\)'' assigns the specific value \(2\) (newly introduced) to the previously defined general variable \(c\). Conversely, ``\(c := 2\)'' would mean introducing a new variable \(c\) and explicitly defining it to be \(2\).}
\[
a<b \xRightarrow[\text{with } c=:2]{P} 2a<2b
\]
may be regarded as irreducible.\footnote{The term ``base field'' is inspired by polynomial irreducibility, which also depends on the choice of base field.} In a more advanced base field, we may directly admit the specialized shortcut trail
\[
\times 2:\quad x<y\Rightarrow 2x<2y
\]
as atomic, so that
\[
a<b \xRightarrow{\times 2} 2a<2b
\]
is irreducible as well. Since the advanced base field is larger, it naturally contains at least all trails from the elementary base field, and may also include additional ones; in particular, it is natural to regard the shortcut trail \(\times 2\) as atomic. This is formalized in Lemma~\ref{lem:B-extension-irreducible}.\qed
\end{example}

\begin{remark}\label{rem:trail-general-route-unit}
A \textit{trail} may itself be viewed as a (general) \textit{route unit}. A concrete route unit is then an instance obtained by renaming symbols and instantiating parameters. For example, the trail \(P\) above can be written as
\[
P:\ (x<y,\ c>0)\xRightarrow{P} cx<cy,
\]
and
\[
a<b \xRightarrow[\text{with }x=:a,\ y=:b,\ c=:2]{P} 2a<2b
\]
is an instance of it.
\end{remark}

The phrase ``allowing only parameter-substitution/renaming'' in Definition~\ref{def:B-irreducible} is formalized by the following lemma.

\begin{lemma}[Instantiation lemma]\label{lem:instantiation}
Let
\[
A(\mathbf{u}) \xRightarrow{P(\mathbf{u})} B(\mathbf{u})
\]
be a parameterized route unit whose validity is guaranteed whenever a hypothesis family \(H(\mathbf{u})\) holds.
If \(\mathbf{u}=: \mathbf{a}\) and \(H(\mathbf{a})\) holds, then
\[
A(\mathbf{a}) \xRightarrow{P(\mathbf{a})} B(\mathbf{a})
\]
is a valid route unit.
\end{lemma}

\begin{proof}
\[
\boxed{A(\mathbf{u}) \xRightarrow{P(\mathbf{u})} B(\mathbf{u}) \text{ whenever } H(\mathbf{u}) \text{ holds}}
\xRightarrow{\mathbf{u}=: \mathbf{a}}
\boxed{A(\mathbf{a}) \xRightarrow{P(\mathbf{a})} B(\mathbf{a}).}
\]
\end{proof}

The next lemma says that once a route unit is irreducible, enlarging the base field does not make it reducible.

\begin{lemma}[Monotonicity under base-field extension]\label{lem:B-extension-irreducible}
If \(\mathbb{B}\subseteq \mathbb{B}'\) and a route unit \(A\xRightarrow{P}B\) is \(\mathbb{B}\)-irreducible, then it is also \(\mathbb{B}'\)-irreducible.
\end{lemma}

\begin{proof}
\[
\begin{aligned}
\boxed{A\xRightarrow{P}B} \text{ is }\mathbb{B}\text{-irreducible}
&\xLeftrightarrow{\text{def}}
P\in\textcolor{blue}{\mathbb{B}}\text{ (and the single-application condition holds)}\\
&\xRightarrow{\textcolor{blue}{\mathbb{B}}\subseteq\mathbb{B}'\text{ given}}
P\in\mathbb{B}'\text{ (and the single-application condition holds)}\\
&\xLeftrightarrow{\text{def}}
\boxed{A\xRightarrow{P}B} \text{ is }\mathbb{B}'\text{-irreducible}. 
\end{aligned}\qedhere
\]
\end{proof}

In Example \ref{ex:ineq-basefield}, if we denote the elementary base field by \(\mathbb{B}\) and the advanced base field by \(\mathbb{B}'\), then \(\mathbb{B}\subseteq \mathbb{B}'\) and
\[
a<b \xRightarrow[\text{with }x=:a,\ y=:b,\ c=:2]{P} 2a<2b
\]
is \(\mathbb{B}\)-irreducible, hence also \(\mathbb{B}'\)-irreducible by the lemma. By contrast,
\[
a<b \xRightarrow{\times 2} 2a<2b
\]
may be \(\mathbb{B}'\)-irreducible without being \(\mathbb{B}\)-irreducible.

Choosing an appropriate pedagogical base field is crucial in mathematical teaching. Since it is impractical to reduce every step to fundamental definitions, the base field should typically include the relevant prerequisites together with the material already established for the intended audience. For example, when teaching algebraic geometry, a suitable base field at a given stage might include the necessary background from abstract algebra and commutative algebra, together with the algebraic geometry already developed up to that point. Most importantly, even if a trail is treated as atomic in the base field, the route unit should not be left defective: the trail should still be explicitly stated.

\subsection{Routeways}\label{Sec:routeways}

\begin{definition}[Routeway]\label{def:routeway}
A \textbf{routeway} (or simply \textbf{route}) is a finite sequence of route units
\[
A_0 \xRightarrow{P_1} A_1 \xRightarrow{P_2} \cdots \xRightarrow{P_n} A_n
\]
in which the terminal waypoint of each route unit is the initial waypoint of the next. Such a routeway may be denoted by
\[
A_0 \rightsquigarrow A_n.
\]
The waypoint \(A_0\) is the \textbf{starting point}, and \(A_n\) is the \textbf{destination}.

A routeway is \textbf{\(\mathbb{B}\)-irreducible} if each of its route units is \(\mathbb{B}\)-irreducible.
\end{definition}

\begin{corollary}\label{cor:B-extension-routeway}
If \(\mathbb{B}\subseteq \mathbb{B}'\), then every \(\mathbb{B}\)-irreducible routeway is \(\mathbb{B}'\)-irreducible.
\end{corollary}

\begin{proof}
\[
\begin{aligned}
\boxed{\gamma:
A_0 \xRightarrow{P_1} A_1 \xRightarrow{P_2} \cdots \xRightarrow{P_n} A_n}
\text{ is }\mathbb{B}\text{-irreducible}
&\xRightarrow{\text{Definition~\ref{def:routeway}}}
\boxed{A_{i-1}\xRightarrow{P_i}A_i} \text{ is }\mathbb{B}\text{-irreducible for each }i\\
&\xRightarrow{\text{Lemma~\ref{lem:B-extension-irreducible}}}
\boxed{A_{i-1}\xRightarrow{P_i}A_i} \text{ is }\mathbb{B}'\text{-irreducible for each }i\\
&\xRightarrow{\text{Definition~\ref{def:routeway}}}
\gamma \text{ is }\mathbb{B}'\text{-irreducible}.
\end{aligned}
\]
\end{proof}

Equivalently, once we form the route graph \(\Gamma_{\mathbb{B}}\), a \(\mathbb{B}\)-irreducible routeway is exactly a directed walk (edge-path) in \(\Gamma_{\mathbb{B}}\).

\begin{definition}[Route graph]\label{def:routegraph}
Fix a pedagogical base field \(\mathbb{B}\). The \textbf{route graph} \(\Gamma_{\mathbb{B}}\) is the labeled directed multigraph (quiver) whose vertices are waypoints and whose directed edges are the \(\mathbb{B}\)-irreducible route units
\[
A \xRightarrow{P} B
\qquad (P\in\mathbb{B}).
\]
Different trails from the same initial waypoint \(A\) to the same terminal waypoint \(B\) are regarded as distinct edges. A two-way route unit
\[
A \xLeftrightarrow{P} B
\]
contributes the two directed edges
\[
A \xRightarrow{P} B
\qquad\text{and}\qquad
B \xRightarrow{P} A.
\]
\end{definition}

By convention, for every waypoint \(A\), we allow the empty routeway \(A\rightsquigarrow A\), containing no route units.

We prefer $\mathbb{B}$-irreducible routeways because, relative to a fixed base field \(\mathbb{B}\), they require no extra cognitive effort beyond applying atomic trails in \(\mathbb{B}\).

\begin{example}\label{e1}
An algebra professor might state:\footnote{In this example, we assume $p$ is a prime number.}
\begin{quote}
``Since the order of the group \( G \) is \( p^2 \), it is solvable.''
\end{quote}

This sentence corresponds to a single-route-unit routeway:
\[
|G| = p^2 \Longrightarrow G \text{ is solvable,}
\]
which is clearly \textit{defective}, as it lacks the explicit trail (justification). 

For the professor, recognizing gaps in their own teaching can be difficult, as the logical structure has become fully internalized; they no longer need to consciously identify the trail to know where it lies. For an expert auditing the lecture, it’s easy to follow along, as they instinctively supply the missing logical steps in their minds. But for a beginner unfamiliar with group theory, the sentence is likely to appear entirely opaque and meaningless.

Suppose instead the professor states:
\begin{quote}
“We know the order of the group $G$ is $p^2$. Since every abelian group is solvable, it follows that $G$ is solvable.”
\end{quote}
First, it’s unclear what the \textit{initial waypoint} is and where the \textit{trail} begins. Second, this essentially corresponds to:
\[
    |G| \text{ (the order of } G) = p^2 
    \xRightarrow{\text{Every abelian group is solvable}}  
    G \text{ is solvable.}
\]
Although it is not \textit{defective} anymore, it is \textit{reducible}, as this justification does not directly state, “Given the \textit{initial waypoint} ($|G| = p^2$), we have the \textit{terminal waypoint} ($G \text{ is solvable}$).” 

Instead, an \textit{irreducible routeway} should be explicitly structured as:

\begin{tcolorbox}[width=\linewidth, colframe=black, colback=white, rounded corners, boxrule=1pt]
\[
\begin{aligned}
|G|=p^2
&\overset{P_1}{\Longrightarrow}
G \cong \mathbb{Z}_{p^2}\ \text{or}\ \mathbb{Z}_p \times \mathbb{Z}_p\\
&\overset{P_2}{\Longrightarrow}
G \text{ is abelian}\\
&\overset{P_3}{\Longrightarrow}
G \text{ is solvable.}
\end{aligned}
\]

\vspace{0.5em}
\textbf{Where:}
\begin{itemize}
\item[$P_1$:] Every group of order \(p^2\) is isomorphic to \(\mathbb{Z}_{p^2}\) or \(\mathbb{Z}_p \times \mathbb{Z}_p\).
\item[$P_2$:] Both \(\mathbb{Z}_{p^2}\) and \(\mathbb{Z}_p \times \mathbb{Z}_p\) are abelian; hence any group isomorphic to either is abelian.
\item[$P_3$:] Every abelian group is solvable.
\end{itemize}
\end{tcolorbox}

We can still describe the routeway in words, but we should make sure no trails are left out. For example, to describe \[
A \ (\textbf{Initial Waypoint}) \ \xRightarrow{P \ (\textbf{Trail})} \ B \ (\textbf{Terminal Waypoint}),
\]
verbally, we can say ``Since $A$, we have $B$. This is because of $P$.'' Yet omitting even one trail forces students to spend a tremendous amount of time reconstructing it on their own. Worse still, many lectures—and a significant portion of textbooks and papers—omit nearly every trail, rendering most route units defective by default.

With these three connected route units powered by $P_1$, $P_2$, and $P_3$ (properties proved earlier in the lecture, making the route units non-defective), students won’t be stuck asking, “How did you get that?” They’ll see every step. Even if they don't know what abelian groups and solvable groups are at all, they can easily follow this irreducible route unit:
\[
G \text{ is abelian}  
    \xRightarrow{\text{Every abelian group is solvable}}  
    G \text{ is solvable.}
\]
\qed
\end{example}

\begin{example}[Finite integral domain \(\Rightarrow\) field]\label{e1a}

When proving that every finite integral domain is a field, a professor might write on the board:\footnote{Here, an integral domain means a commutative ring with \(1\neq 0\) and no zero divisors.}
\begin{quote}
Let \(R\) be a finite integral domain and \(0\neq a\in R\). Define \(m_a:R\to R\) by \(m_a(x)=ax\). Then \(m_a\) is injective. Since \(R\) is finite, \(m_a\) is surjective. Thus \(1\in \operatorname{Im}(m_a)\), so \(\exists b\in R\) with \(m_a(b)=1\), i.e.\ \(ab=1\). Since \(R\) is commutative, also \(ba=1\), so \(b\) is an inverse of \(a\). Therefore \(R\) is a field.
\end{quote}

This is correct, but as before, it compresses several trails into a few words. Below we rewrite it as an explicit routeway with each trail labeled.

As above, define \(m_a:R\to R\) by \(m_a(x)=ax\). Then

\begin{tcolorbox}[width=\linewidth, colframe=black, colback=white, rounded corners, boxrule=1pt]
\[
\begin{aligned}
&a\in R,\ a\neq 0,\ R \text{ is a finite integral domain}\\
&\xRightarrow{P_1} \forall x,y\in R,\ ax=ay \Rightarrow x=y\\
&\xLeftrightarrow{\text{def of }m_a} \forall x,y\in R,\ m_a(x)=m_a(y)\Rightarrow x=y\\
&\xLeftrightarrow{\text{def of injective}} m_a \text{ is injective}\\
&\xLeftrightarrow[R\text{ is finite (given)}]{P_2} m_a \text{ is surjective}\\
&\xRightarrow[\text{with }y=:1]{P_3 (\text{def of surjective})} \exists b\in R:\ m_a(b)=1\\
&\xLeftrightarrow{\text{def of }m_a}
\exists b\in R:\ 
\underset{\mathclap{\hspace{7.5em}\begin{array}[t]{@{}c@{\;}l@{}}
\rotatebox[origin=c]{90}{$=$} & \text{\scriptsize commutativity of }R\\[-0.4ex]
ba &
\end{array}}}{ab}=1.
\end{aligned}
\]
\vspace{0.5em}
\textbf{Where:}
\begin{itemize}
\item[$P_1$:] In an integral domain, every \(a\neq 0\) is not a zero divisor (equivalently, \(ax=ay\Rightarrow x=y\)).
\item[$P_2$:] On a finite set, a map is injective if and only if it is surjective.
\item[$P_3$:] If \(f:X\to Y\) is surjective, then for every \(y\in Y\) there exists \(x\in X\) such that \(f(x)=y.\)
\item[$P_4$:] A commutative ring with \(1\) in which every nonzero element has a multiplicative inverse is a field.
\end{itemize}
\end{tcolorbox} 
Since \(a\neq 0\) was arbitrary, the box above gives
\[
\forall a\in R,\ a\neq 0 \Rightarrow \exists b\in R:\ ab=ba=1.
\]
Hence \(P_4\) implies that \(R\) is a field. \qed
\end{example}

In irreducible routeways, every logical connection is explicitly revealed. Nothing remains hidden or ambiguous. Such clarity profoundly transforms mathematical teaching and learning, empowering both teachers and students. Explicitly structured routeways should replace mostly defective descriptive narratives in textbooks and lectures, fundamentally changing mathematics education for the better.

\section{Familiar Anchors: Connecting to Known Concepts}\label{Sec:anchors}

New concepts should never be introduced in isolation. An effective way to understand unfamiliar concepts is to anchor them to familiar ideas.

\begin{definition}[Anchor]\label{def:anchor}
An \textbf{anchor} is a designated waypoint, fixed in advance for the intended audience, used as a reference point for reaching new waypoints.
\end{definition}

When convenient, we may also speak informally of \emph{anchor trails}, meaning familiar trails already available in the pedagogical base field. The formal distance and closure constructions below, however, use only anchor waypoints.

This prompts us to expand our scope: not merely concentrating on individual route units, but on various approaches and alternative routeways. In this section, rather than examining individual route units separately, we treat the entire routeway as a cohesive unit and broaden our perspective accordingly.

\begin{definition}[Roadmap]\label{def:roadmap}
Fix \(\Gamma_{\mathbb{B}}\). A \textbf{roadmap} is a chosen collection of routeways in \(\Gamma_{\mathbb{B}}\) from a given starting point to a given destination.
\end{definition}

\begin{remark}[Roadmaps as subgraphs]
Fix \(\mathbb{B}\). A roadmap can be viewed as the subgraph of \(\Gamma_{\mathbb{B}}\) (Definition~\ref{def:routegraph}) consisting of all waypoints and route units that appear in the collection of known routeways.
\end{remark}

For instance, the following figure illustrates a roadmap comprising two distinct routeways.

\begin{figure}[H]
    \centering
     \includegraphics[width=0.7\textwidth]{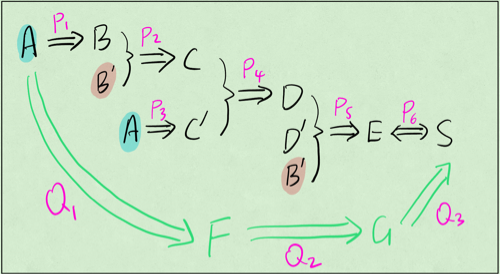}
    \caption{Roadmap\protect\footnotemark}
    \label{fig:roadmap-example}
\end{figure}
\footnotetext{The four colors depicted in Figure~\ref{fig:roadmap-example} are officially named: {\color[HTML]{14C7DE}aqua} (blue), {\color[HTML]{E26A6A}bittersweet shimmer} (red), {\color[HTML]{2FDA77}mermaid's kiss} (green, indicating a new routeway), and {\color[HTML]{FF00D4}razzle dazzle rose} (pink, for trails \(P_i\) and \(Q_i\)). Using these whimsical color names can significantly enhance teaching engagement.}

\begin{definition}[Teaching atlas]\label{def:atlas}
Fix \(\Gamma_{\mathbb{B}}\) and a set \(X\) of target waypoints. A \textbf{teaching atlas} (or simply \textbf{atlas}) on \(X\) is a collection of roadmaps
\[
\mathcal{A}=\{\mathcal{R}_i\}_{i\in I}
\]
in \(\Gamma_{\mathbb{B}}\), each possibly with a different starting point and destination, such that for every \(x\in X\), there exist \(i\in I\) and a routeway in \(\mathcal{R}_i\) in which \(x\) appears as a waypoint.
\end{definition}

Thus an atlas is a collection of roadmaps whose combined coverage reaches the target waypoints.

Our focus now shifts from simply ensuring each routeway is clear, to determining the \textit{optimal routeway} to utilize during instruction. The objective is efficiency: classroom time is limited, and exploring every conceivable routeway is impractical. Selecting a routeway anchored to what students already know significantly improves comprehension and retention.  This mirrors Polya’s problem-solving philosophy~\cite{polya}, where success often hinges on connecting the unfamiliar to the familiar.

Consider the process of language learning. When an English speaker inquires, “What does ‘gato’ mean in Spanish?”, the most confusing response would be: “It's a small domesticated feline known for its agility and independence. It can typically climb trees and move swiftly along narrow structures...”. In contrast, the optimal response is simply: “It means ‘cat’.” Immediate clarity is achieved without extraneous details; understanding follows naturally from familiarity. 
There is no need anymore to hear the whole part about climbing trees; once we know it means ``cat,'' we know it can climb trees.

\begin{figure}[H]
    \centering
    \includegraphics[width=\textwidth]{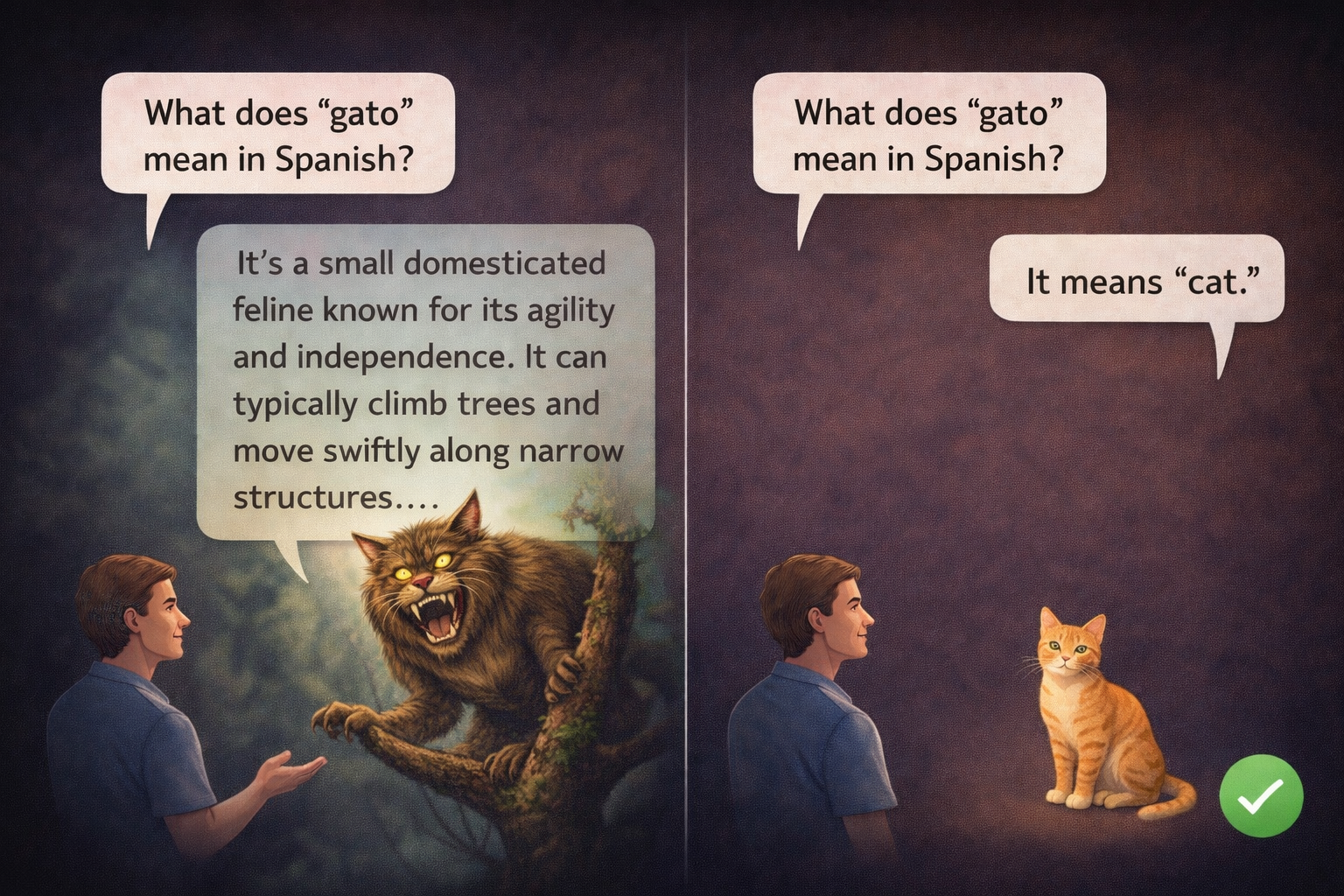}
    \caption{Descriptive overload versus immediate anchoring to a familiar concept.}
    \label{fig:gato}
\end{figure}

Nevertheless, educators frequently adopt the overly complicated first approach. Rather than saying, “A \textit{module} is just like a vector space, but over a ring instead of a field,'' they often begin with an intricate abstract definition, forcing students to grapple with unfamiliarity. Similarly, instead of succinctly stating, “Take a module, add multiplication like in a ring, and you get an \textit{algebra},” educators might immediately delve into formal definitions involving bilinear maps and associativity, thereby unnecessarily complicating the initial introduction of the concept. For students already familiar with modules and rings, to achieve a quick grasp, the first version of the definition of algebras should be stated first while the formal details can come later.
\begin{example}\label{e2}
Imagine teaching a child how to get to school. The child already knows two pieces of information clearly:
\begin{enumerate}
    \item How to travel from home to their Father’s Office.
    \item How to proceed from their Father’s Office to the Campus Gate.
\end{enumerate}

Now you must instruct the child on the complete route to the school. The most confusing way is to provide overly detailed, step-by-step instructions such as:
\begin{quote}
``Walk ten feet, turn left, go two miles, turn right, then take another left after the third intersection. Or maybe after the fifth—depends on traffic. Oh, and if that road is blocked, you might need to backtrack. Actually, there’s another way too, but it’s a bit longer. So it's better to take a cab. Have you taken one before? Well, let me teach you how to do it…''
\end{quote}
Not only the child, but even adults who regularly commute to the school and know the route thoroughly would find this explanation challenging to follow. The optimal explanation would be:
\begin{quote}
``First, go to your dad’s office. From there, go to the campus gate, just like before. Then walk straight until the end of the road, turn right, and then you’ll see your school.”
\end{quote}
Immediate clarity is provided, without cognitive overload. This approach is not only straightforward to follow but also easy to remember. The previously described \textit{most confusing way} resembles an inefficient professor attempting to explain the longer routeway involving the trails denoted by \(P\) in Figure \ref{fig:roadmap-example}, while the \textit{optimal way} corresponds to utilizing the shorter (``{\color[HTML]{2FDA77}mermaid's kiss}'') routeway with trails labeled \(Q\), connecting familiar anchors {\color[HTML]{2FDA77}$F$} ({\color[HTML]{2FDA77}F}ather’s Office), {\color[HTML]{2FDA77}$G$} (Campus {\color[HTML]{2FDA77}G}ate), and ultimately arriving at {\color[HTML]{2FDA77}$S$} ({\color[HTML]{2FDA77}S}chool).\footnote{Interestingly, these letters naturally match the initial letters of the corresponding locations.}

More importantly, once the child reaches their destination, a shift occurs: the paths become internalized. The unfamiliar transforms into something familiar. The School is no longer a mystery—it becomes a new anchor, the next ``Father’s Office.'' From that point forward, the child will no longer merely follow directions blindly; they will perceive alternate routes and discover shortcuts independently. \qed
\end{example}

\begin{example}\label{e2'}
In a probability class, a professor is teaching this problem:\footnote{This problem is adapted from Exercise 3.10.1 in \cite{grimmett}.}
\begin{quote}
    Consider a symmetric simple random walk \( S \) with \( S_0 = 0 \). Recall that
\[
T_x := \inf \{ n \geq 1 : S_n = x \} \text{ (where $S_n$ is the position of the random walk at step $n$)}
\]
is the \textit{first hitting time} of state \( x \). Let \( T := T_0 \) be the first return time to the origin. Prove, for \(m\ge 1\):
\begin{equation}
     \mathbb{P}(T = 2m) = \frac{1}{2m - 1} \binom{2m}{m} \left(\frac{1}{2}\right)^{2m}.
    \label{eq:main-prob}
\end{equation}
\textbf{Hint:} 
\begin{equation}\label{eq:hint-prob}
\mathbb{P}(T = 2m)
=
\mathbb{P}\Big(S_{2m-1} = 1,\ S_k \neq 0 \ \forall\, 1 \le k < 2m\Big).
\tag{\makebox[0pt][r]{(}\theequation\makebox[0pt][l]{)$'$}}
\end{equation}
That is, the probability we want is the same as
\[
\mathbb{P}(S_{2m-1} = 1 , \, \text{no touching 0 before step } 2m \text{ after leaving } (0,0)),
\]
because to return to 0 for the first time at step \(2m\), the walk must be at \(+1\) or \(-1\) at step \(2m-1\) and must never touch 0 before; by symmetry the \(+1\) and \(-1\) cases contribute equally, and the final step from \(\pm1\) to \(0\) contributes a factor \(\tfrac12\). We may take this hint as given in this problem.
\end{quote}

Assume students have already learned and proved key theorems from prior lectures—our familiar anchors.\footnote{To connect with common textbook forms, explanations are included below, but for this example, we can focus only on the boxed equations \eqref{eq:ballot-count} and \eqref{eq:path-count-general} and ignore their deductions. That said, their deductions provide good examples of how we write irreducible routeways.}

\vspace{0.5cm}
\paragraph{\textbf{Anchor 1: The Ballot Problem (Father's Office)}}

Suppose candidate \(A\) receives \(\alpha\) votes and candidate \(B\) receives \(\beta\) votes, where \(\alpha>\beta\). The probability that \(A\) always stays ahead of \(B\) during the counting process is given by \eqref{eq:bal-geo}:

\stepcounter{footnote}\edef\fnA{\number\value{footnote}}
\stepcounter{footnote}\edef\fnB{\number\value{footnote}}

\begin{equation}
\hspace*{-5.0em}
\begin{array}{cc}
    \mathbb{P}(\text{A always ahead}) = \dfrac{\alpha - \beta}{\alpha + \beta} 
= 
  \mathbb{P}\left(\text{path from } (0,0) \to (\alpha + \beta, \alpha - \beta) \text{ without hitting \(0\) after time \(0\)}\ \middle|\ S_{\alpha+\beta}=\alpha-\beta\right)\footnotemark[\fnA]
\\
    \begin{array}{c}
   \hspace{4.5cm}\raisebox{0.6ex}{$\Vert$\scriptsize{\text{i.e.}}}\\
  \hspace{4.5cm}\dfrac{\#(\text{paths from } (0,0) \to (\alpha + \beta, \alpha - \beta) \text{ without hitting \(0\) after time \(0\)})}{\#(\text{paths from } (0,0) \to (\alpha + \beta, \alpha - \beta))}
\end{array}
\end{array}
\label{eq:bal-geo}
\end{equation}
\footnotetext[\fnA]{The latter equation is the geometric interpretation of the probability in the Ballot problem, easily arising from its proof using the idea of random walks.}

\noindent\hspace*{-3.0em}
\makebox[\textwidth][l]{%
$\displaystyle
\begin{rcases*}
\eqref{eq:bal-geo} \xLeftrightarrow{\text{i.e.}}
\dfrac{\alpha - \beta}{\alpha + \beta} =
\dfrac{\#(\text{paths from } (0,0) \to (\alpha + \beta, \alpha - \beta) \text{ without hitting \(0\) after time \(0\)})}{\#(\text{paths from } (0,0) \to (\alpha + \beta, \alpha - \beta))}
\\
\text{Denominator\footnotemark[\fnB]:  } \#(\text{paths from } (0,0) \to (\alpha + \beta, \alpha - \beta)) \xlongequal[a=:0,\ n=:\alpha+\beta,\ b=:\alpha-\beta]{\eqref{eq:path-binomial}\;\text{ below with}}
\displaystyle{\binom{\alpha + \beta}{\frac{(\alpha + \beta) + (\alpha - \beta)}{2}}}
\end{rcases*}\xRightarrow[\text{by Denominator}]{\text{multiply \eqref{eq:bal-geo}}}
$
}
\footnotetext[\fnB]{In the following equation, it might be tempting to simply write `` $\xlongequal[]{\eqref{eq:path-binomial}\;\text{ below}}$'' here (skipping the substitution details below the long equals sign). But this omission can slow the reader down. They may need to pause, compare symbols, and match this step with the general fact (trail), which is not ideal. To fully qualify the definition of an \textit{irreducible route unit}, every step should be direct, with no unnecessary pauses or extra effort.}

\[
\text{Numerator:  }\#(\text{paths from } (0,0) \to (\alpha + \beta, \alpha - \beta) \text{ without hitting \(0\) after time \(0\)}) =\dfrac{\alpha - \beta}{\alpha + \beta} \binom{\alpha + \beta}{\frac{(\alpha + \beta )+ (\alpha - \beta)}{2}}
\]

\begin{equation}
    \hspace*{-3.0em} \xRightarrow{\text{let }\alpha+\beta=:p,\ \alpha-\beta=:q}
\boxed{\#(\text{paths from } (0,0) \to (p,q) \text{ without hitting \(0\) after time \(0\)}) = \frac{q}{p} \binom{p}{\frac{p+q}{2}}\qquad (p\ge 1,\ q>0)}
\tag{\makebox[0pt][r]{(}\theequation\makebox[0pt][l]{)$'$}}\quad
\label{eq:ballot-count}
\end{equation}
Note that we adopt the convention
\[
\binom{n}{r}:=0
\qquad\text{if } r\notin\mathbb{Z}\text{ or } r<0 \text{ or } r>n.
\]

\vspace{0.5cm}

\paragraph{\textbf{Anchor 2: Generalized Path Counting Formula (Campus Gate)}}

The probability of a symmetric random walk reaching \( (n, b) \) from \( (0, a) \), possibly under additional constraints (such as not touching zero), is given by:
\begin{equation}
\boxed{
  \mathbb{P}_a\big(S_n = b,\ \text{constraints}\big)
  =
  \#(\text{valid paths from } (0, a) \to (n, b)) \left(\frac{1}{2}\right)^n}
\label{eq:path-count-general}
\end{equation}

Note that we usually omit the subscript $a$ in $  \mathbb{P}_a$ if $a=0$.

In particular, for unrestricted paths:
\begin{equation}
  \mathbb{P}_a(S_n = b) = \#(\text{paths from } (0, a) \to (n, b)) \times \left(\frac{1}{2}\right)^n
\label{eq:path-count-unrestricted}
\tag{\makebox[0pt][r]{(}\theequation\makebox[0pt][l]{)$'$}}
\end{equation}
where
\begin{equation}
\#(\text{paths from } (0,a) \to (n,b)) =
\begin{cases}
\displaystyle \binom{n}{\frac{n+b-a}{2}}, & n+b-a \text{ is even},\\[0.6em]
0, & n+b-a \text{ is odd}.
\end{cases}
\label{eq:path-binomial}
\tag{\makebox[0pt][r]{(}\theequation\makebox[0pt][l]{)$''$}}
\end{equation}
\vspace{0.5cm}
\paragraph{\textbf{Without Anchors}}  

The professor might explain the problem like this:

\begin{quote}
\textbf{Professor:} ``Let's just follow the hint and analyze this step by step. We need to count the number of paths reaching \( S_{2m-1} = 1 \) while avoiding zero before step \( 2m-1 \). 

First, to avoid hitting $0$, the very first step is forced to be $+1$. As a result, the walk begins at $(1,1)$ rather than $(0,0)$. Consider all possible random walk paths from \( (1,1) \) to \( (2m-1,1) \). The horizontal movement is $(2m-1)-1=2m-2$ steps. Since each step is either \( +1 \) or \( -1 \), the total number of such paths is given by the binomial coefficient \( \binom{2m-2}{m-1} \). Do you follow? Well, it’s just because there are \( 2m-2 \) total steps, and \( \frac{2m-2}{2}=m-1 \) of them must be \( +1 \) steps, the rest \( -1 \). So, this binomial coefficient is the total number of these paths.

But we only want the paths that never touch zero after leaving the start. To count how many violate this condition, we use the reflection principle: every path that touches zero corresponds uniquely to a path that, instead of ending at \( (2m-1,1) \), ends at \( (2m-1,-1) \). The number of such invalid paths is \( \binom{2m-2}{m-2} \) since now there must be $m$ down $-1$ steps and $m-2$ up $+1$ steps. So, subtracting the invalid cases, the number of valid paths is \( \binom{2m-2}{m-1} - \binom{2m-2}{m-2} \).

Now, let’s talk probability. Each individual step in the random walk has a probability of \( \frac{1}{2} \). That means any specific sequence of steps—up, down, up, up, down, whatever—has probability \( 2^{-(2m-2)} \). This holds for every possible path, valid or not. So, each valid path occurs with probability \( 2^{-(2m-2)} \). That makes sense, right? Any sequence of \( 2m-2 \) steps is just a random string of up and down moves, and each move is equally likely. But wait—our first step still has to be \( +1 \) to reach \( S_{1} = 1 \), which happens exactly half the time. You see where that comes from? We already fixed \( S_{1} = 1 \), which forced our first step, so we just tack on the extra factor of \( \frac{1}{2} \). That means we need to multiply by \( \frac{1}{2}\), giving us \( 2^{-2m+1} \). 

Now we multiply by the number of valid paths to get the final probability
\[
  \mathbb{P}(T = 2m) = \left(\binom{2m-2}{m-1} - \binom{2m-2}{m-2}\right) 2^{-2m+1}.
\]
Finally, using the binomial identity \( \binom{2m-2}{m-1} - \binom{2m-2}{m-2} = \frac{1}{2(2m-1)} \binom{2m}{m} \), which you can easily verify, we arrive at
\[
  \mathbb{P}(T = 2m) = \frac{1}{2m-1} \binom{2m}{m} 2^{-2m}.
\]
And that’s our result.''

\textbf{Student:} ``Is the proof like the Ballot theorem we learned? Can we use the Ballot theorem to solve it?''

\textbf{Professor:} ``Yes, the idea is similar because they both use the reflection principle. But they're not the same. In this problem, we are not dealing with votes but with first return times. The Ballot theorem is about staying ahead, while here we’re controlling the first time we return to the starting point. But in a sense, both of them can be thought of as a random walk without crossing $0$. So, they’re related. Does that answer your question?''

\textbf{Student:} ``Yes, thanks.''
\end{quote}
This may sound formal and feel detailed—but it is overwhelming, like reinventing the wheel to prove something already known in a slightly different way. The process is akin to guiding the child to their Father's Office using an unfamiliar, unnecessarily complex route in Example \ref{e2}, instead of simply saying, “Go to your father's office.” Thus, even if the professor’s route units are mostly \textit{irreducible}, this is not the best routeway in the roadmap: it fails to make use of familiar anchors. Even if it is easy to follow, it remains hard for students to remember or reconstruct on their own.

\vspace{0.5cm}
\paragraph{\textbf{With Anchors}}  
\begin{enumerate}
    \item \textbf{Home → Father’s Office}: Using \eqref{eq:ballot-count} with \(p=:2m-1,\ q=:1\)    \[
    \#(\text{valid paths}) = \frac{1}{2m - 1} \binom{2m-1}{m}.
    \]

    \item \textbf{Father’s Office → Campus Gate}: Using \eqref{eq:path-count-general} with \(a=:0,\ n=:2m-1,\ b=:1\), $\text{constraints}=:S_k \neq 0 \ \forall \ 1 \le k < 2m$
 \begin{equation}\label{eq:path-count-special}
\mathbb{P}\Big(S_{2m-1} = 1,\ S_k \neq 0 \ \forall\, 1 \le k < 2m\Big)
=
\frac{1}{2m - 1} \binom{2m-1}{m} \left(\frac{1}{2}\right)^{2m-1}.
\end{equation}

    \item \textbf{Campus Gate → School}:  
\[\hspace*{-5.0em}
\begin{array}{cc}
 \eqref{eq:path-count-special} \xLeftrightarrow{\text{i.e.}}\displaystyle \mathbb{P}(S_{2m-1} = 1, \,  S_k \neq 0 \ \forall \ 1 \le k < 2m) = \frac{1}{2m - 1} \textcolor{magenta}{\binom{2m-1}{m}}\left(\frac{1}{2}\right)^{2m-1}
\\ \begin{array}{c}
   \hspace{5.0cm}\raisebox{1.5ex}{$\Vert$\scriptsize{\eqref{eq:hint-prob}}}\\
      \hspace{5.5cm} \displaystyle\mathbb{P}(T=2m)\quad\xlongequal[]{\smiley }
\end{array}
    \begin{array}{c}
\hspace{0.9cm}\raisebox{0.1ex}{$\Vert$ {\scriptsize \text{binomial identity:} $\displaystyle{\textcolor{magenta}{\binom{2m-1}{m}}=\frac{1}{2}\binom{2m}{m} }$}}\\
  \hspace{-3.0cm}\quad\displaystyle\frac{1}{2m - 1} \binom{2m}{m}\left(\frac{1}{2}\right)^{2m}
\end{array}
\end{array}
\]
\end{enumerate}
\qed
\end{example}

That illustrates how math should be taught with efficiency: Start with what students already know, and then build on it. New knowledge isn’t an island; it’s an extension of something familiar.

\section{Mathematical Compass and Driving Simulation: Trailblazing New Routeways}

With the insights gained from the previous sections, we have become educators capable of making mathematics highly accessible. However, we should aim further. This section explores tools that can permanently enhance our students' mathematical thinking and their academic landscape.

A roadmap tells us where to go, but what if the road doesn't exist yet? Venturing into uncharted territory requires tools for discovery and exploration. In unexplored mathematical landscapes, we need a \textit{compass} (motivation) for direction and a \textit{driving simulation} (examples) to experience the route before traveling independently.

\subsection{Mathematical Compass (Motivation)}

Let's go back to the route unit
\[
A \ (\textbf{Initial Waypoint}) \ \xRightarrow{P \ (\textbf{Trail})} \ B \ (\textbf{Terminal Waypoint}).
\]
With \( P \) explicitly stated, we have essentially answered the question “How do we get \( B \) from \( A \) in this step?” for a broad range of students. However, even with the trail made clear, a thoughtful student might still ask, “Why use \( P \) here? What led to that choice?” That is a great sign—it means they’re not lost in \textit{how}, and have moved on to \textit{why}—what motivates us to think of choosing this trail here. With such motivation, even in the absence of explicit guidance, students will have both the handrail and the confidence to navigate the route on their own. If we, as educators, can go a step further—motivating each transition with clarity and concrete examples, even before students ask—we transcend competent instruction. We become exceptional educators who cultivate mathematical intuition and inspire confident exploration.

As Schoenfeld argues in~\cite{schoenfeld}, making motivations (the ``why'' part) clear means that we not only lead students to the destination clearly (\S \ref{Sec:trails} Clear Trails) and effectively (\S \ref{Sec:anchors} Familiar Anchors), but also give them the necessary tools—a \textit{compass}—to recover when they go off track, discover new waypoints, and expand the roadmap.

\begin{definition}[Mathematical Compass]
A \textbf{mathematical compass} (or simply \textbf{compass}) is the underlying motivation guiding the choice of a particular route unit within a routeway. 
\end{definition}

\begin{remark}
A compass is meta-mathematical information attached to a route unit or routeway which explains why that step or route is a natural choice. Unlike trails, route distance, or closure, a compass is not part of the formal graph-theoretic structure; it shows motivation rather than validity.
\end{remark}

Once students know the rough direction of the destination (which they usually do, since the problem goal or expectation is often given or known), they can navigate on their own—exploring and completing routeways built by others, creating new routeways, and eventually arriving at their destination.

For example, in Example~\ref{e2}, when telling the child to go to their Father’s Office, a curious child might ask why go there. We can tell them: ``Because your dad’s office is a familiar place to you, just like your home, and you know how to get to many places from there. So, starting from a familiar place makes things easier. That's why you can first go there.'' This illustrates \textit{why} we take those steps—not \textit{how} (which is already known to the child), but \textit{why}, in the big picture of getting to School. The same applies to Example~\ref{e2'}.

\begin{example}
In general,
\[
\lim\limits_{n \to \infty} a_n = L
\xLeftrightarrow{\text{def}}
(\forall \varepsilon>0)(\exists N\in\mathbb{N})(\forall n\ge N)\,|a_n-L|<\varepsilon.
\]
To use the definition of limits to prove
\[
\lim\limits_{n \to \infty} q^n = 0, \quad \text{for } 0 < |q| < 1,
\]
we can teach it like this: Given \(\varepsilon>0\), choose an integer \(N\) such that\footnote{It exists by the Archimedean property.}
\[
N > \frac{\log \varepsilon}{\log |q|}.
\]
Then, since \(\log|q|<0\), for every \(n\ge N\) we have
\[
n > \frac{\log \varepsilon}{\log |q|}
\Longrightarrow
n\log|q|<\log\varepsilon
\Longrightarrow
|q|^n<\varepsilon.
\]
Hence \(|q^n-0|<\varepsilon\).
Thus, the proof follows directly from the definition.

Notice this explanation corresponds to the following irreducible\footnote{In a standard calculus or analysis base field, this is irreducible.} route unit:
\[
(\forall \varepsilon>0)(\exists N\in\mathbb{N})(\forall n\ge N)\,|q^n-0|<\varepsilon
\]
\[
\xLeftrightarrow[\text{with } L =: 0,\ a_n =: q^n]{\text{def}}
\lim\limits_{n \to \infty} q^n = 0.
\]
That's clear enough. Nevertheless, the threshold \( \frac{\log \varepsilon}{\log |q|} \) may seem to appear out of nowhere. Why do we compare \(N\) to this seemingly arbitrary quantity? Even though there is no logical gap in the proof, explaining the motivation behind choosing \( N \) helps students develop intuition:

\begin{quote}
To ensure \( |q^n - 0| < \varepsilon \), we solve for \( n \):
\[
|q^n - 0| < \varepsilon \quad \Leftrightarrow \quad n > \frac{\log \varepsilon}{\log |q|}.
\]
This inspires us to choose \(N\) to be any integer strictly larger than \(\frac{\log \varepsilon}{\log |q|}\), ensuring that for any \(n \ge N\), we get \( |q^n - 0| < \varepsilon \).
\end{quote}

Providing this insight encourages students to develop similar reasoning skills, helping them not only understand the content but also gain a deeper sense of mathematical thinking. \qed
\end{example}

\subsection{Driving Simulation (Concrete Examples)}

Concrete and fundamental examples, serving as special cases of general mathematical facts, are essential in math teaching. This parallels the idea of allowing students to explore routes through a driving simulation, becoming familiar with the terrain and landmarks before navigating high-stakes real-world situations. When students ultimately traverse these paths independently, they can compare their experience with prior simulations, enriching their understanding significantly.

A driving simulation informally means previewing or validating a general routeway through concrete special cases (examples). We now formalize this idea.

\begin{definition}[Parameterized routeway and driving simulation]\label{def:simulation}
A \textbf{parameterized routeway} is a routeway whose waypoints and trails depend on formal parameters \(\mathbf{u}\).

A \textbf{driving simulation} is a concrete specialization
\[
\mathbf{u}=: \mathbf{a}
\]
of such a routeway, provided the instantiated hypotheses of all trails are satisfied.
\end{definition}

\begin{corollary}[Simulation preserves routeways]\label{lem:simulation-preserves}
Let \(\gamma(\mathbf{u})\) be a parameterized routeway valid whenever \(H(\mathbf{u})\) holds.
If \(\mathbf{u}=: \mathbf{a}\) and \(H(\mathbf{a})\) holds, then the specialized routeway \(\gamma(\mathbf{a})\) is valid.
\end{corollary}

\begin{proof}
Each route unit of \(\gamma(\mathbf{u})\) specializes to a valid route unit by Lemma~\ref{lem:instantiation}. Hence the specialized route units concatenate to a valid routeway \(\gamma(\mathbf{a})\).
\end{proof}

Conversely, an invalid simulation provides a rigorous test for refuting a proposed generalized statement.

\begin{corollary}[Concrete counterexample detector]\label{cor:counterexample-detector}
If a proposed universal trail or routeway fails under a specialization
\[
\mathbf{u}=: \mathbf{a}
\]
that satisfies its stated hypotheses, then the proposed universal statement is false.
\end{corollary}

\begin{proof}
For a single route unit, this is the contrapositive of Lemma~\ref{lem:instantiation}.
For a whole routeway, it is the contrapositive of Corollary~\ref{lem:simulation-preserves}.
\end{proof}

Mazur eloquently asks “When is one thing equal to some other thing?”~\cite{mazur}, illustrating how names, notation, and abstraction often hide rather than reveal meaning. Abstraction serves to unify and generalize: for instance, \( 2 + 1 \) may denote both ``2 cats and 1 dog'' or ``2 roses and 1 dog''—its meaning abstracted from the objects themselves. However, abstraction without a solid grounding resembles a floating tree crown: disconnected, lacking roots or a trunk, and thus incapable of growth. In reality, abstraction is always supported by concreteness. Concrete examples anchor and nourish abstractions, allowing them to flourish meaningfully.  This disconnect between structure and meaning also reflects concerns raised in~\cite{thom} about modern mathematics becoming detached from human experience.

\begin{figure}[H]
    \centering
    \includegraphics[width=\textwidth]{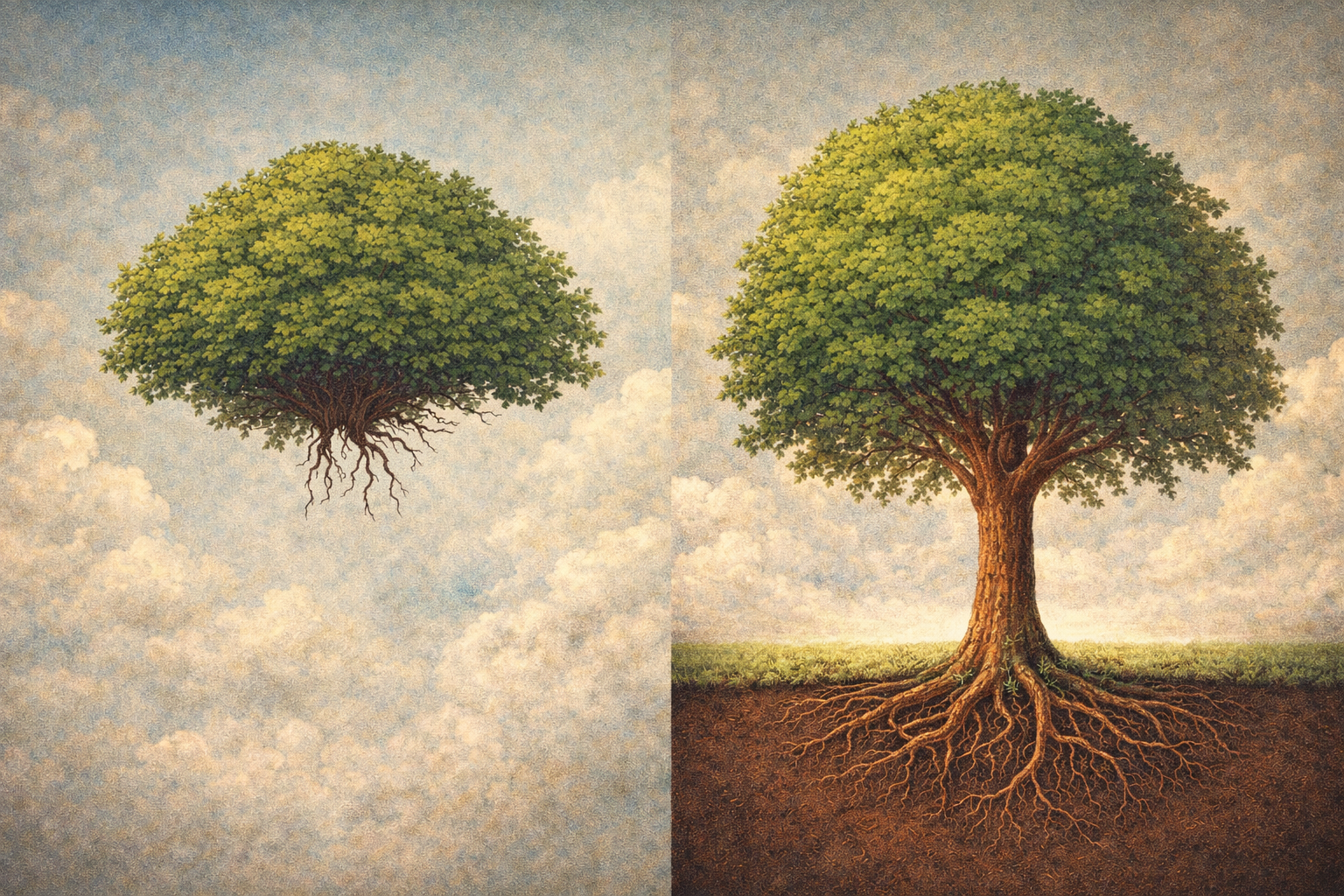}
    \caption{Abstraction without grounding versus abstraction rooted in concrete examples.}
    \label{fig:roots}
\end{figure}

For example, the definition of a topology is not arbitrary—it generalizes fundamental properties of open sets in metric spaces \cite{munkres}. The structure is as follows:

\begin{quote}
\textbf{Definition.} A \textit{topology} on a set \( X \) is a collection \( \mathcal{T} \) of subsets of \( X \) satisfying the following properties:
\begin{enumerate}
    \item \( \varnothing \in \mathcal{T} \) and \( X \in \mathcal{T} \),
    \item The union of the elements of any subcollection of \( \mathcal{T} \) is in \( \mathcal{T} \),
    \item The intersection of the elements of any finite subcollection of \( \mathcal{T} \) is in \( \mathcal{T} \).
\end{enumerate}
A set \( X \) for which a topology \( \mathcal{T} \) has been specified is called a \textit{topological space}. If \( X \) is a topological space with topology \( \mathcal{T} \), we say that a subset \( U \subseteq X \) is an \textit{open set} if \( U \in \mathcal{T} \). Using this terminology, a topological space is a set \( X \) together with a collection of subsets of \( X \), called open sets, such that \( \varnothing \) and \( X \) are open, and such that arbitrary unions and finite intersections of open sets are open.
\end{quote}

These conditions did not arise arbitrarily; they were abstracted from the behavior of open sets in metric spaces. In any metric space, the collection of open sets satisfies these three properties above: containing the empty set and whole space, being closed under arbitrary unions, and being closed under finite intersections. These properties were then taken as axioms to generalize the notion of ``open sets'' beyond metric spaces, forming the foundation of topology. Without understanding this origin, the abstract definition appears detached and loses intuitive significance, as words floating without a root.

Likewise, instead of first introducing the orbit in group theory as
\[
G \cdot x := \{ g \cdot x \mid g \in G \}
\]
with examples, one could begin directly with a categorical formulation:

\begin{quote}
The orbit of \( s \in S \) is the image of the morphism \( o(\cdot, s): G \to S \), induced by the action map \( a: G \times S \to S \) with \( g \mapsto a(g, s) \).
\end{quote}

While mathematically precise, such a definition offers little intuition. It does not explain why the term “orbit” is used, nor does it convey any concretely motivating idea. Such abstraction is valuable—indeed, this is the spirit in which category theory arose, unifying disparate branches of mathematics—but it must be grounded in a deep understanding of the intuitive special cases it abstracts from.

\begin{example}
To illustrate function parity, consider the textbook definition of an even function, stated rigorously as:
\begin{quote}
A function \( f \) is \textit{even} if and only if \( f(-x) = f(x) \) for all \( x \) in the domain of \( f \).
\end{quote}
This definition may make perfect sense to someone already familiar with the concept of function parity, but for a beginner, the rigorous form alone conveys little. A student encountering the definition \( f(-x) = f(x) \) for the first time is unlikely to grasp its meaning without seeing it in action:
\begin{quote}
Let's consider \( f(x) = x^2 \). Let \( x = 3 \). Then \( f(3) = 3^2 = 9 \). Now take \( x = -3 \), and observe that \( f(-3) = (-3)^2 = 9 \) as well. That is, \( f(3) = f(-3) \). In fact, we can easily verify that this holds for all \( x \in \mathbb{R} \): \( f(-x) = f(x) \). Such behavior is striking and consistent, so it deserves a name.

More generally, for any power function of the form \( f(x) = x^n \), if \( n \) is an even integer (e.g., \( -2, 0, 4 \)), then \( f(-x) = f(x) \) wherever it is defined. This observation motivates the terminology: any function—whether a power function or not—that satisfies \( f(-x) = f(x) \) is called an \emph{even function}. 
\end{quote}
Thus, these abstract definitions and names are not arbitrary; they are rooted in concrete examples and patterns, giving meaning and structure to the general concept. \qed
\end{example}

In mathematics education, abstraction should always be grounded in concrete examples. Without them, abstract concepts remain isolated and forgettable. With them, abstraction acquires meaning, intuition, and permanence in learners' minds.

\section{Route Geometry}\label{Sec:route-geometry}

Once routeways are viewed as directed paths in \(\Gamma_{\mathbb{B}}\), they admit a natural geometry. This section defines route length, refinement, directed distance, geodesic routeways, anchor distance, route intervals, anchor excess, and closure, and studies how the geometric notions depend on the pedagogical base field.

\subsection{Geodesic Structure}\label{Sec:geodesic-structure}

Once several routeways are available, it is natural to compare them quantitatively rather than only descriptively.

\begin{definition}[Route length, route distance, and geodesic routeway]\label{def:geodesic}
The \textbf{length} of a routeway \(\gamma\) is the number of route units in \(\gamma\), denoted \(\ell(\gamma)\). That is, for a routeway
\[
\gamma:
A_0 \xRightarrow{P_1} A_1 \xRightarrow{P_2} \cdots \xRightarrow{P_n} A_n,
\]
we have
\[
\ell(\gamma):=n.
\]

Fix a pedagogical base field \(\mathbb{B}\). For waypoints \(A,B\), define the \textbf{route distance}\footnote{In general, \(d_{\mathbb{B}}(A,B)\neq d_{\mathbb{B}}(B,A)\). Thus \(d_{\mathbb{B}}\) is a directed distance rather than a metric in the usual symmetric sense.}

\[
d_{\mathbb{B}}(A,B):=
\begin{cases}
\displaystyle \min_{\gamma:\,A\rightsquigarrow B\text{ in }\Gamma_{\mathbb{B}}}\ell(\gamma), & \text{if such a routeway exists},\\[0.6em]
+\infty, & \text{otherwise.}
\end{cases}
\]

A \textbf{geodesic routeway} from \(A\) to \(B\) is a routeway \(\gamma^\ast:A\rightsquigarrow B\) with
\[
\ell(\gamma^\ast)=d_{\mathbb{B}}(A,B)
\]
whenever \(d_{\mathbb{B}}(A,B)<+\infty\).
\end{definition}

\begin{lemma}[Existence of geodesics]\label{lem:geodesic-existence}
If \(d_{\mathbb{B}}(A,B)<+\infty\), then there exists a geodesic routeway from \(A\) to \(B\).
\end{lemma}

\begin{proof}
\[
\begin{aligned}
d_{\mathbb{B}}(A,B)<+\infty
&\xRightarrow{\text{def}}
\varnothing\neq \{\ell(\gamma): \gamma:A\rightsquigarrow B \text{ in }\Gamma_{\mathbb{B}}\}\subseteq \mathbb{Z}_{\ge 0}\\
&\xRightarrow{\text{well-ordering}}
\exists\,\gamma^*:A\rightsquigarrow B:\ \ell(\gamma^*)=d_{\mathbb{B}}(A,B),
\end{aligned}
\]
where the well-ordering principle says that every nonempty subset of \(\mathbb{Z}_{\ge 0}\) has a least element. Hence \(\gamma^*\) is geodesic by definition.
\end{proof}
Note that, by Definition~\ref{def:routegraph}, every routeway in \(\Gamma_{\mathbb{B}}\) is automatically \(\mathbb{B}\)-irreducible.

\begin{lemma}[Additivity of route length]\label{lem:length-additive}
If
\[
\gamma_1:A\rightsquigarrow B
\qquad\text{and}\qquad
\gamma_2:B\rightsquigarrow C
\]
are routeways, then
\[
\ell(\gamma_2\circ\gamma_1)=\ell(\gamma_1)+\ell(\gamma_2).
\]
\end{lemma}

\begin{proof}
Write
\[
\gamma_1:
A_0 \xRightarrow{P_1} A_1 \xRightarrow{P_2} \cdots \xRightarrow{P_m} A_m,
\]
\[
\gamma_2:
B_0 \xRightarrow{Q_1} B_1 \xRightarrow{Q_2} \cdots \xRightarrow{Q_n} B_n,
\]
with \(A_m=B_0\). Then
\[
\gamma_2\circ\gamma_1:
A_0 \xRightarrow{P_1} A_1 \xRightarrow{P_2} \cdots \xRightarrow{P_m} A_m=B_0 \xRightarrow{Q_1} B_1 \xRightarrow{Q_2} \cdots \xRightarrow{Q_n} B_n,
\]
so
\[
\ell(\gamma_2\circ\gamma_1)
\xlongequal{\text{def}}
m+n
\xlongequal{\text{def}}
\ell(\gamma_1)+\ell(\gamma_2).
\]
\end{proof}

\subsection{Refinement}\label{Sec:refinement}

A central claim of this paper is that narrative proofs should be replaced by explicit routeways. The next definitions formalize that replacement process.

\begin{definition}[Refinement of routeways]\label{def:refinement}
Let \(\gamma\) and \(\eta\) be finite routeways with the same starting point and destination.
We say that \(\eta\) is a \textbf{refinement} of \(\gamma\), and write
\[
\gamma \preceq \eta,
\]
if \(\eta\) is obtained from \(\gamma\) by replacing zero or more route units, each with endpoints \(X\) and \(Y\), by nonempty finite routeways from \(X\) to \(Y\).
\end{definition}

\begin{remark}
If \(\gamma \preceq \eta\), then \(\ell(\gamma)\le \ell(\eta)\). Thus refinement preserves or increases the number of route units; it never shortens a routeway.
\end{remark}

\begin{lemma}[Refinement is a preorder]\label{lem:refinement-preorder}
The relation \(\preceq\) is reflexive and transitive.
\end{lemma}

\begin{proof}
Reflexivity:
\[
\hspace*{-5.0em}
\gamma\preceq\gamma
\xLeftrightarrow{\text{def}}
\gamma \text{ is obtained from }\gamma\text{ by replacing route units by nonempty finite routeways with the same endpoints,}
\]
which holds by replacing no route units.

Transitivity:
\[
\begin{aligned}
\gamma\preceq\eta \text{ and } \eta\preceq\theta
&\xLeftrightarrow{\text{def}}
\begin{cases}
\eta \text{ is obtained from }\gamma\text{ by local routeway replacements}\\
\theta \text{ is obtained from }\eta\text{ by local routeway replacements}
\end{cases}\\
&\xRightarrow{\text{compose}}
\theta \text{ is obtained from }\gamma\text{ by local routeway replacements}\\
&\xLeftrightarrow{\text{def}}
\gamma\preceq\theta.
\end{aligned}
\]
\end{proof}

\noindent This shows \(\preceq\) is a preorder on finite routeways with fixed starting point and destination.

\begin{example}[Refinement is not symmetric]
Assume the displayed route units are valid, and let
\[
\gamma:\quad A \xRightarrow{P} B,
\qquad
\eta:\quad A \xRightarrow{P_1} C \xRightarrow{P_2} B.
\]
Then
\[
\gamma \preceq \eta,
\]
since the single route unit from \(A\) to \(B\) may be replaced by the two-step routeway from \(A\) to \(B\) through \(C\). However,
\[
\eta \npreceq \gamma,
\]
because every refinement of \(\eta\) still contains \(C\) as a waypoint: the first replacement block ends at \(C\), and the second begins at \(C\). Hence no refinement of \(\eta\) can equal the one-step routeway \(\gamma\). Thus refinement is not symmetric.
\end{example}

\begin{definition}[Presentation-equivalence]\label{def:presentation-equivalence}
For finite routeways \(\gamma\) and \(\eta\) with the same starting point and destination, define
\[
\gamma\sim\eta
\quad\Longleftrightarrow\quad
\gamma\preceq\eta \text{ and } \eta\preceq\gamma.
\]
In this case, we say that \(\gamma\) and \(\eta\) are \textbf{presentation-equivalent}.
\end{definition}

By Lemma~\ref{lem:refinement-preorder}, together with symmetry from the definition, the relation \(\sim\) is an equivalence relation.

\begin{remark}[Pedagogical meaning of presentation-equivalence]
Presentation-equivalence means that \(\gamma\) and \(\eta\) differ only by local choices of presentation: each can be obtained from the other by replacing some route units with routeways having the same endpoints. In particular, it does not identify a single route unit with a longer routeway through new intermediate waypoints.
\end{remark}

Crucially, if every individual step of a routeway can be refined into irreducible units, the entire routeway admits an irreducible refinement.

\begin{lemma}[Concatenation preserves irreducibility]\label{lem:concat-irreducible}
Let
\[
\eta=\eta_1\eta_2\cdots\eta_n
\]
be a routeway. If each \(\eta_i\) is \(\mathbb{B}\)-irreducible, then \(\eta\) is \(\mathbb{B}\)-irreducible.
\end{lemma}

\begin{proof}
\[
\begin{aligned}
\eta_i \text{ is }\mathbb{B}\text{-irreducible for each }i
&\xRightarrow{\text{Definition~\ref{def:routeway}}}
\text{every route unit of each }\eta_i\text{ is }\mathbb{B}\text{-irreducible}\\
&\xRightarrow{\eta=\eta_1\eta_2\cdots\eta_n}
\text{every route unit of }\eta\text{ is }\mathbb{B}\text{-irreducible}\\
&\xRightarrow{\text{Definition~\ref{def:routeway}}}
\eta \text{ is }\mathbb{B}\text{-irreducible}.
\end{aligned}
\]
\end{proof}

\begin{theorem}[Existence of irreducible refinement]\label{thm:irreducible-refinement}
Fix a pedagogical base field \(\mathbb{B}\). Let
\[
\gamma=e_1e_2\cdots e_n
\]
be a finite routeway, where each \(e_i\) is a valid route unit, viewed as a length-one routeway. Assume that for each \(i\), there exists a finite \(\mathbb{B}\)-irreducible routeway \(\eta_i\) such that
\[
e_i\preceq \eta_i.
\]
Then \(\gamma\) admits a finite \(\mathbb{B}\)-irreducible refinement.
\end{theorem}

\begin{proof}
Write
\[
\gamma=e_1e_2\cdots e_n,
\qquad
e_i\preceq \eta_i
\quad(1\le i\le n),
\]
where each \(\eta_i\) is finite and \(\mathbb{B}\)-irreducible. Since \(e_i\preceq \eta_i\), by Definition \ref{def:refinement}, each \(\eta_i\) has the same starting point and destination as \(e_i\). 
Therefore, by definition
\[
\eta:=\eta_1\eta_2\cdots \eta_n
\]
is a finite routeway, and
\[
\gamma=e_1e_2\cdots e_n \preceq \eta.
\]
By Lemma~\ref{lem:concat-irreducible}, \(\eta\) is \(\mathbb{B}\)-irreducible. Therefore \(\eta\) is a finite \(\mathbb{B}\)-irreducible refinement of \(\gamma\).
\end{proof}

\subsection{Anchors and Closure}\label{Sec:anchor-geometry}

\begin{definition}[Anchor distance]\label{def:anchor-distance}
Fix a nonempty set \(S\) of waypoints serving as anchors. For a waypoint \(x\), define
\[
d_{\mathbb{B}}(S,x):=\min_{s\in S} d_{\mathbb{B}}(s,x),
\]
where the minimum is taken in \(\mathbb{Z}_{\ge 0}\cup\{+\infty\}\).
\end{definition}

\begin{lemma}[More anchors never increase anchor distance]\label{lem:anchor-enlarge}
If \(S\subseteq T\), then
\[
d_{\mathbb{B}}(T,x)\le d_{\mathbb{B}}(S,x)
\]
for every waypoint \(x\).
\end{lemma}

\begin{proof}
\[
\begin{aligned}
S\subseteq T
&\xRightarrow{\text{def}}
\{d_{\mathbb{B}}(t,x):t\in T\}\supseteq \{d_{\mathbb{B}}(s,x):s\in S\}\\
&\xRightarrow{\min \text{ in }[0,+\infty]}
\min_{t\in T}d_{\mathbb{B}}(t,x)\le \min_{s\in S}d_{\mathbb{B}}(s,x)\\
&\xLeftrightarrow{\text{def}}
d_{\mathbb{B}}(T,x)\le d_{\mathbb{B}}(S,x).
\end{aligned}
\]
\end{proof}

For any pedagogical base field \(\mathbb{F}\), if
\[
\gamma: A\rightsquigarrow B
\]
is a routeway in \(\Gamma_{\mathbb{F}}\), then by the definition of \(d_{\mathbb{F}}(A,B)\),
\begin{equation}\label{eq:dandl}
d_{\mathbb{F}}(A,B)\le \ell(\gamma).
\end{equation}

\begin{lemma}[Base field extension shrinks route distance]\label{lem:extension}
If \(\mathbb{B}\subseteq \mathbb{B}'\), then
\[
d_{\mathbb{B}'}(A,B)\le d_{\mathbb{B}}(A,B)
\]
for all waypoints \(A,B\).
\end{lemma}

\begin{proof}
If \(d_{\mathbb{B}}(A,B)=+\infty\), then
\[
d_{\mathbb{B}'}(A,B)\le d_{\mathbb{B}}(A,B)
\]
is immediate.

Otherwise, Lemma~\ref{lem:geodesic-existence} provides a geodesic routeway \(\gamma\) from \(A\) to \(B\) in \(\Gamma_{\mathbb{B}}\). Then
\[
\begin{aligned}
\gamma: A\rightsquigarrow B \text{ in }\Gamma_{\mathbb{B}}
&\xRightarrow{\text{Corollary~\ref{cor:B-extension-routeway}}}
\gamma: A\rightsquigarrow B \text{ in }\Gamma_{\mathbb{B}'}\\
&\xRightarrow[\text{with } \mathbb{F}=:\mathbb{B}']{\eqref{eq:dandl}}
d_{\mathbb{B}'}(A,B)\le \ell(\gamma)\xlongequal[\gamma \text{ geodesic in }\Gamma_{\mathbb{B}}]{\text{def}}d_{\mathbb{B}}(A,B)
\end{aligned}
\]
\end{proof}

This monotonicity suggests modeling learning as an expanding pedagogical base field.

\begin{definition}[Knowledge filtration]
A \textbf{knowledge filtration} is an increasing sequence
\[
\mathbb{B}_0 \subseteq \mathbb{B}_1 \subseteq \cdots \subseteq \mathbb{B}_t \subseteq \cdots
\]
where \(\mathbb{B}_t\) is the pedagogical base field at time \(t\) (e.g.\ week \(t\)).
\end{definition}

\begin{remark}
As \(t\) increases, route distance \(d_{\mathbb{B}_t}(A,B)\) is nonincreasing (Lemma~\ref{lem:extension}).
\end{remark}

Route distances naturally satisfy the triangle inequality, which is crucial for estimating the length of a detour via an anchor.

\begin{lemma}[Triangle inequality]\label{lem:triangle}
For any waypoints \(A,B,C\),
\[
d_{\mathbb{B}}(A,C)\le d_{\mathbb{B}}(A,B)+d_{\mathbb{B}}(B,C).
\]
\end{lemma}

\begin{proof}
If \(d_{\mathbb{B}}(A,B)=+\infty\) or \(d_{\mathbb{B}}(B,C)=+\infty\), then
\[
d_{\mathbb{B}}(A,C)\le d_{\mathbb{B}}(A,B)+d_{\mathbb{B}}(B,C)
\]
is immediate.

Otherwise, by Lemma~\ref{lem:geodesic-existence}, choose geodesic routeways
\[
\gamma_1^*:A\rightsquigarrow B,
\qquad
\gamma_2^*:B\rightsquigarrow C
\]
in \(\Gamma_{\mathbb{B}}\). Then
\[
\begin{aligned}
\boxed{\gamma_1^*:A\rightsquigarrow B,\quad \gamma_2^*:B\rightsquigarrow C}
&\xRightarrow{\text{concatenate}}
\boxed{\gamma_2^*\circ\gamma_1^*:A\rightsquigarrow C}\\
&\xRightarrow{\eqref{eq:dandl}}
d_{\mathbb{B}}(A,C)\le \ell(\gamma_2^*\circ\gamma_1^*)\\
&\hphantom{d_{\mathbb{B}}(A,C)}\xlongequal{\text{Lemma~\ref{lem:length-additive}}}
\ell(\gamma_1^*)+\ell(\gamma_2^*)\\
&\hphantom{d_{\mathbb{B}}(A,C)}\xlongequal[\gamma_1^*,\,\gamma_2^*\text{ geodesic}]{\text{def}}
d_{\mathbb{B}}(A,B)+d_{\mathbb{B}}(B,C).
\end{aligned}
\]
\end{proof}

A second basic property of route distance is that geodesicity is inherited by subroutes.

\begin{lemma}[Optimal substructure of geodesics]\label{thm:subroute}
Let \(\gamma:A\rightsquigarrow C\) be a geodesic routeway. If
\[
\gamma=\gamma_2\circ\gamma_1,
\qquad
\gamma_1:A\rightsquigarrow B,
\qquad
\gamma_2:B\rightsquigarrow C,
\]
then both \(\gamma_1\) and \(\gamma_2\) are geodesic.
\end{lemma}

\begin{proof}
Assume
\[
\gamma=\gamma_2\circ\gamma_1,
\qquad
\gamma_1:A\rightsquigarrow B,
\qquad
\gamma_2:B\rightsquigarrow C.
\]

If \(\gamma_1\) were not geodesic, then, by definition,
\begin{equation}\label{eq:exsubroute}
\exists\,\gamma_1':A\rightsquigarrow B
\quad\text{with}\quad
\ell(\gamma_1')<\ell(\gamma_1).
\end{equation}
Hence
\[
\begin{aligned}
\gamma_1':A\rightsquigarrow B,\ \gamma_2:B\rightsquigarrow C
&\xRightarrow{\text{concatenate}}
\gamma_2\circ\gamma_1':A\rightsquigarrow C\\
&\xRightarrow[\text{take }\mathbb{F}=\mathbb{B}]{\eqref{eq:dandl}}
d_{\mathbb{B}}(A,C)
\le \ell(\gamma_2\circ\gamma_1')\\
&\hphantom{d_{\mathbb{B}}(A,C)\le{}}\xlongequal{\text{Lemma~\ref{lem:length-additive}}}
\ell(\gamma_1')+\ell(\gamma_2)\\
&\hphantom{d_{\mathbb{B}}(A,C)\le{}}
\overset{\eqref{eq:exsubroute}}{<}\ell(\gamma_1)+\ell(\gamma_2)\\
&\hphantom{d_{\mathbb{B}}(A,C)\le{}}\xlongequal[\gamma=\gamma_2\circ\gamma_1]{\text{Lemma~\ref{lem:length-additive}}}
\ell(\gamma)\\
&\hphantom{d_{\mathbb{B}}(A,C)\le{}}\xlongequal[\gamma\text{ geodesic}]{\text{def}}
d_{\mathbb{B}}(A,C),
\end{aligned}
\]
a contradiction. Thus \(\gamma_1\) is geodesic.

Similarly, \(\gamma_2\) is geodesic.
\end{proof}

This leads to a precise notion of when a waypoint lies on a shortest route from \(A\) to \(B\).

\begin{definition}[Route interval]\label{def:route-interval}
Assume \(d_{\mathbb{B}}(A,B)<+\infty\). The \textbf{route interval} from \(A\) to \(B\) is
\[
I_{\mathbb{B}}(A,B)
:=
\left\{
F:
d_{\mathbb{B}}(A,B)=d_{\mathbb{B}}(A,F)+d_{\mathbb{B}}(F,B)
\right\}.
\]
\end{definition}

The next theorem identifies route intervals exactly as the waypoints lying on a geodesic routeway from \(A\) to \(B\).

\begin{theorem}[Interval characterization of perfect anchors]\label{thm:route-interval}
Assume \(d_{\mathbb{B}}(A,B)<+\infty\). For a waypoint \(F\), the following are equivalent:
\begin{enumerate}
\item \(F\in I_{\mathbb{B}}(A,B)\).
\item There exists a geodesic routeway from \(A\) to \(B\) passing through \(F\).
\end{enumerate}
\end{theorem}

\begin{proof}
\((1)\Rightarrow(2)\):  
\[
F\in I_{\mathbb{B}}(A,B)
\xRightarrow{\text{def}}
d_{\mathbb{B}}(A,F)<+\infty \text{ and } d_{\mathbb{B}}(F,B)<+\infty.
\]
By Lemma~\ref{lem:geodesic-existence}, choose geodesic routeways
\[
\gamma_1^*:A\rightsquigarrow F,
\qquad
\gamma_2^*:F\rightsquigarrow B.
\]
Then
\[
\begin{aligned}
F\in I_{\mathbb{B}}(A,B)
\xLeftrightarrow{\text{def}}
d_{\mathbb{B}}(A,B)
&=d_{\mathbb{B}}(A,F)+d_{\mathbb{B}}(F,B)\\
&\xlongequal[\gamma_1^*,\,\gamma_2^*\text{ geodesic}]{\text{def}}
\ell(\gamma_1^*)+\ell(\gamma_2^*)\\
&\xlongequal{\text{Lemma~\ref{lem:length-additive}}}
\ell(\gamma_2^*\circ\gamma_1^*).
\end{aligned}
\]
Hence \(\gamma_2^*\circ\gamma_1^*\) is a geodesic routeway from \(A\) to \(B\) passing through \(F\).

\((2)\Rightarrow(1)\):
Suppose there exists a geodesic routeway \(\gamma^*:A\rightsquigarrow B\) passing through \(F\). Write
\[
\gamma^*=\gamma_2\circ\gamma_1,
\qquad
\gamma_1:A\rightsquigarrow F,
\qquad
\gamma_2:F\rightsquigarrow B.
\]
Thus, by Lemma~\ref{lem:length-additive}, 
\begin{equation}\label{eq:lroute-interval}
    \ell(\gamma^\ast)=\ell(\gamma_1)+\ell(\gamma_2).
\end{equation}
Then
\[
\begin{aligned}
\gamma^*:A\rightsquigarrow B \text{ geodesic through }F
&\xRightarrow{\text{Lemma~\ref{thm:subroute}}}
\gamma_1:A\rightsquigarrow F,\ \gamma_2:F\rightsquigarrow B \text{ geodesic}\\
&\xRightarrow{\text{def}}
d_{\mathbb{B}}(A,F)=\ell(\gamma_1),\quad d_{\mathbb{B}}(F,B)=\ell(\gamma_2),\quad d_{\mathbb{B}}(A,B)=\ell(\gamma^*)\\
&\xRightarrow{\eqref{eq:lroute-interval}}
d_{\mathbb{B}}(A,B)=d_{\mathbb{B}}(A,F)+d_{\mathbb{B}}(F,B)\\
&\xLeftrightarrow{\text{def}}
F\in I_{\mathbb{B}}(A,B).
\end{aligned}
\]
\end{proof}

The interval condition can also be quantified by measuring the detour created by forcing a route through a chosen waypoint.

\begin{definition}[Anchor excess]\label{def:anchor-excess}
Assume
\[
d_{\mathbb{B}}(A,B),\ d_{\mathbb{B}}(A,F),\ d_{\mathbb{B}}(F,B)<+\infty.
\]
The \textbf{anchor excess} of \(F\) relative to \((A,B)\) is
\[
\operatorname{Ex}_{\mathbb{B}}(F;A,B)
:=
d_{\mathbb{B}}(A,F)+d_{\mathbb{B}}(F,B)-d_{\mathbb{B}}(A,B).
\]
\end{definition}

\noindent By Lemma~\ref{lem:triangle},
\[
\operatorname{Ex}_{\mathbb{B}}(F;A,B)\ge 0.
\]

\begin{corollary}[Zero-excess criterion]\label{cor:anchor-excess-zero}
Assume
\[
d_{\mathbb{B}}(A,B),\ d_{\mathbb{B}}(A,F),\ d_{\mathbb{B}}(F,B)<+\infty.
\]
Then
\[
\operatorname{Ex}_{\mathbb{B}}(F;A,B)=0
\quad\Longleftrightarrow\quad
F\in I_{\mathbb{B}}(A,B).
\]
\end{corollary}

\begin{proof}
\[
\operatorname{Ex}_{\mathbb{B}}(F;A,B)=0
\xLeftrightarrow{\text{def}}
d_{\mathbb{B}}(A,B)=d_{\mathbb{B}}(A,F)+d_{\mathbb{B}}(F,B)
\xLeftrightarrow{\text{def}}
F\in I_{\mathbb{B}}(A,B).
\]
\end{proof}

This suggests two distinguished kinds of waypoints.

\begin{definition}[Perfect and essential waypoints]\label{def:anchor-types}
Assume \(d_{\mathbb{B}}(A,B)<+\infty\). A waypoint \(F\) is called:
\begin{itemize}
\item \textbf{perfect} for \((A,B)\) if \(F\in I_{\mathbb{B}}(A,B)\),
\item \textbf{essential} for \((A,B)\) if every geodesic routeway from \(A\) to \(B\) passes through \(F\).
\end{itemize}
\end{definition}

These notions are related as follows.

\begin{lemma}[Essential implies perfect]\label{lem:anchor-hierarchy}
Assume \(d_{\mathbb{B}}(A,B)<+\infty\). If \(F\) is essential for \((A,B)\), then \(F\) is perfect for \((A,B)\).
\end{lemma}

\begin{proof}
\[F \text{ essential for }(A,B)
\xLeftrightarrow{\text{def}}\text{every geodesic routeway from \(A\) to \(B\) passes through \(F\) } (*)\]
\[
\begin{aligned}
d_{\mathbb{B}}(A,B)<+\infty
&\xRightarrow{\text{Lemma~\ref{lem:geodesic-existence}}}
\exists\text{ geodesic }\gamma^*:A\rightsquigarrow B\\
&\xRightarrow{\text{by }(*)}
\text{\(\gamma^*\) passes through \(F\)}\\
&\xRightarrow{\text{Theorem~\ref{thm:route-interval}}}
F\in I_{\mathbb{B}}(A,B).
\end{aligned}
\]
\end{proof}

\begin{remark}[Why anchors work, formally]
The excess \(\operatorname{Ex}_{\mathbb{B}}(F;A,B)\) measures the detour created by forcing an explanation from \(A\) to \(B\) through \(F\). A perfect waypoint has zero excess, so it lies on a geodesic explanation.
\end{remark}

The next formal ingredient captures the ``once you arrive, it becomes familiar'' phenomenon in Example~\ref{e2}.

\begin{definition}[$\mathbb{B}$-closure of anchors]\label{def:closure}
For a set \(S\) of waypoints serving as anchors, define the \textbf{$\mathbb{B}$-closure}
\[
\mathrm{Cl}_{\mathbb{B}}(S)
:=
\{x:\ \exists s\in S\ \text{and a routeway in }\Gamma_{\mathbb{B}}\text{ from }s\text{ to }x\}.
\]
\end{definition}

\begin{remark}[Interpretation via Example~\ref{e2}]
In Example~\ref{e2}, let
\[
S=\{\text{Father's Office},\ \text{Campus Gate}\}.
\]
Then \(\mathrm{Cl}_{\mathbb{B}}(S)\) is the set of all places the child can reach from these familiar anchors by following routeways in \(\Gamma_{\mathbb{B}}\). For instance, if
\[
\text{Father's Office}\rightsquigarrow \text{Campus Gate}
\qquad\text{and}\qquad
\text{Campus Gate}\rightsquigarrow \text{School},
\]
then
\[
\text{Father's Office}\rightsquigarrow \text{School},
\]
so
\[
\text{School}\in \mathrm{Cl}_{\mathbb{B}}(S).
\]
Thus \(\mathrm{Cl}_{\mathbb{B}}(S)\) includes everything reachable from the familiar anchors.
\end{remark}

\begin{lemma}[Closure axioms]\label{lem:closure}
\(\mathrm{Cl}_{\mathbb{B}}\) is extensive, monotone, and idempotent:
\[
S\subseteq \mathrm{Cl}_{\mathbb{B}}(S),\qquad
S\subseteq T\Rightarrow \mathrm{Cl}_{\mathbb{B}}(S)\subseteq \mathrm{Cl}_{\mathbb{B}}(T),\qquad
\mathrm{Cl}_{\mathbb{B}}(\mathrm{Cl}_{\mathbb{B}}(S))=\mathrm{Cl}_{\mathbb{B}}(S).
\]
\end{lemma}

\begin{proof}
Extensive: each \(s\in S\) reaches itself by the empty routeway.

Monotone: enlarging \(S\) only adds possible starting points.

Idempotent: if \(x\in \mathrm{Cl}_{\mathbb{B}}(\mathrm{Cl}_{\mathbb{B}}(S))\), then \(x\) is reachable from some \(y\) that is reachable from some \(s\in S\); concatenating routeways shows \(x\in \mathrm{Cl}_{\mathbb{B}}(S)\). The reverse inclusion is extensive.
\end{proof}

\section*{Conclusion}

\begin{figure}[H]
    \centering
    \includegraphics[width=1\textwidth]{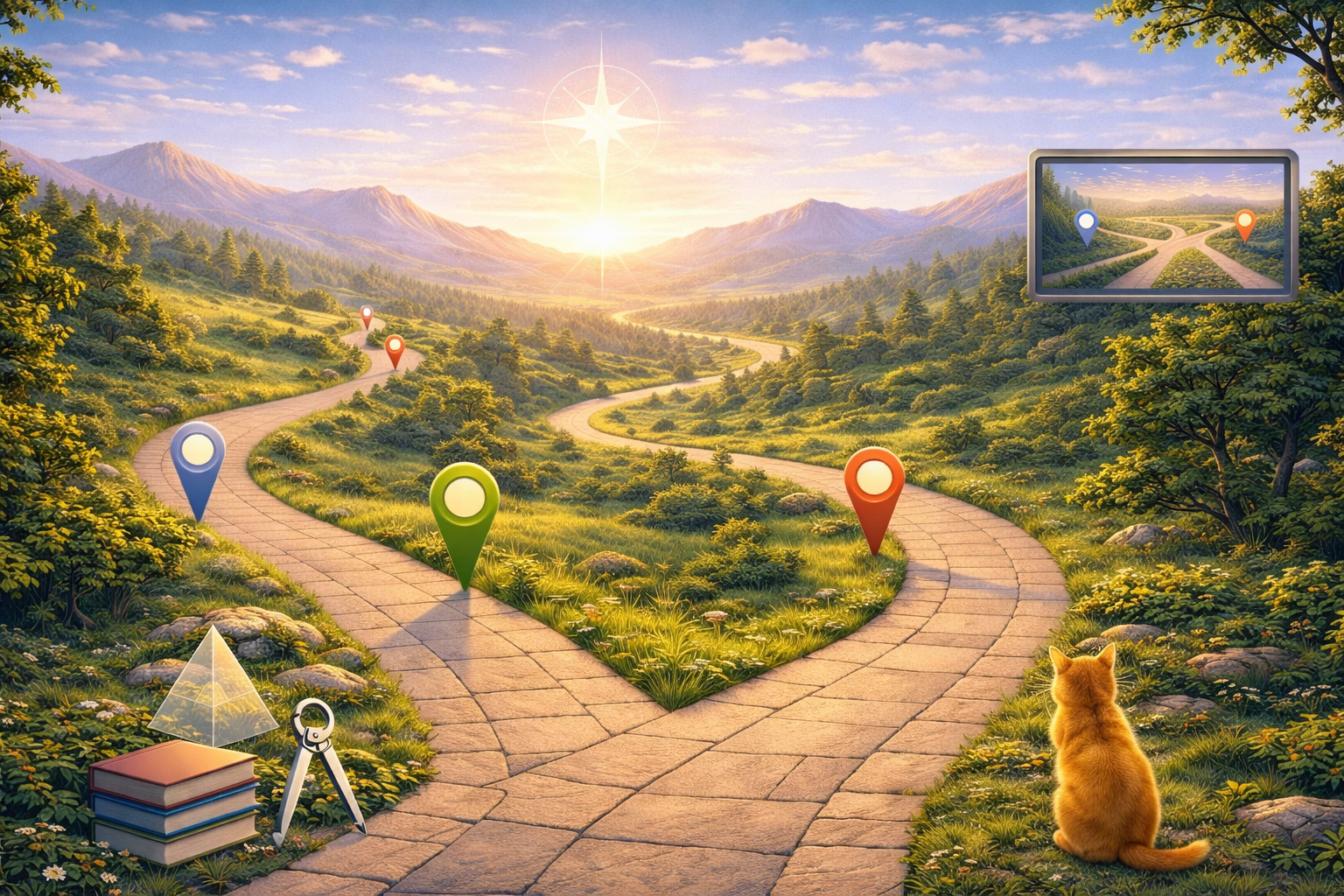}
    \caption{A visual synthesis of routeways, anchors, compass, and driving simulation.}
    \label{fig:atlas-overview}
\end{figure}

The two core principles discussed in the first two sections—\textit{clear trails} and \textit{familiar anchors}—are deeply interconnected. In Example \ref{e1}, we employed three familiar properties \( P_1, P_2, P_3 \) as anchor trails for students. Similarly, in Example \ref{e2}, had we additionally provided an actual roadmap to the child with explicitly annotated paths, the explanation of the way from the Campus Gate to School could have been streamlined significantly: the roadmap itself offering clear guidance. Both principles involve navigating the same mathematical roadmap but emphasize different aspects. Clear trails focus on \textit{visibility}: ensuring each logical step is explicit and clearly presented, leaving no gaps. Familiar anchors focus on \textit{choice}: among various possible paths, those with familiar landmarks should be chosen or prioritized to provide intuitive guidance.

The route geometry developed here has significant potential as a foundational framework. At a formal level, the route graph \(\Gamma_{\mathbb{B}}\) is a labeled directed multigraph, and its routeways form the morphisms of the free category generated by \(\Gamma_{\mathbb{B}}\): waypoints are objects, routeways are morphisms, concatenation is composition, and the empty routeway \(A\rightsquigarrow A\) is the identity at \(A\).

Contrast this with traditional approaches in mathematical instruction:  
Professors often present notes containing largely incomplete route units, selectively filling in missing logical steps during lectures. Lectures thus become compensatory sessions, making up for gaps and ambiguities within the notes. Consequently, textbooks and notes become verbose, structurally ambiguous, and logically opaque.

In response, this paper proposes an alternative instructional framework:

\begin{enumerate}
\item \textbf{Notes should contain explicit, easy-to-follow, irreducible routeways built upon familiar anchors.} If there are any omitted steps in a routeway, they should be intentional, explicitly labeled exercises to encourage active engagement. When structured this way, notes become especially accessible on their own, freeing up valuable lecture time to focus on the following two points.
    \item \textbf{Lectures should emphasize driving simulation.} Rather than solely compensating for incomplete notes, lectures should immerse students through powerful initial experiences with concrete examples or special cases, preparing them for general, abstract concepts.
\item \textbf{Time permitting, each route unit should be guided by a mathematical compass.} By proactively addressing not just \textit{how} but also \textit{why} certain steps are taken, with clear motivations stated, students are equipped with the necessary mathematical compass to navigate independently and creatively through novel problems.

\end{enumerate}

Such an approach makes contemporary mathematics not only accessible but also deeply intuitive. Students thus learn more than mere mathematical techniques; they cultivate critical thinking, creativity, and exploratory skills, expanding the roadmap for future generations of mathematicians.

The current landscape of higher mathematics, in both research and instruction, largely resembles the era of computing before the graphical user interface (GUI). In this analogy, modern mathematical exposition functions like a conceptual DOS: systematic and powerful, yet largely arcane, accessible only to an elite class of initiates. What this paper proposes is not a simplification of mathematical content, but a \textit{GUI for mathematics}: a structural interface that transforms abstract logic into something navigable, teachable, and human. Just as the invention of Windows and Macintosh brought computing into the hands of the many, this framework aims to bring the clarity of advanced mathematics within reach. This is higher math—not for the privileged few, but for the rest of us.

Finally, this paper points out that this conceptual framework extends beyond teaching. In AI-driven mathematical proof analysis and proof generation, clear and explicit roadmaps facilitate logical transparency, verification, and automation. Irreducible routeways simplify complexities, while familiar anchors enhance algorithmic efficiency. Although AI can currently write and analyze mathematical proofs, it frequently errs due to the inherent ambiguity and linguistic flexibility in traditional mathematical writing. Providing structured, rigorously defined route units empowers AI with significantly improved comprehension, generation, and verification of mathematical proofs. Equipped with such structured proof-writing frameworks, algorithms can independently generate mathematically valid proofs, potentially discovering novel routeways previously unexplored by humans through exhaustive verification of logical trails within computationally feasible durations.


\begin{thebibliography}{9}

\bibitem{aristotle-prior}
Aristotle,
\textit{Prior Analytics},
translated by A. J. Jenkinson,
in \textit{The Works of Aristotle}, ed. W. D. Ross,
Oxford University Press, 1928.

\bibitem{lakatos}
I. Lakatos, 
\textit{Proofs and Refutations: The Logic of Mathematical Discovery}, 
Cambridge University Press, 1976.

\bibitem{polya}
G. Pólya, 
\textit{How to Solve It}, 
Princeton University Press, 1945.

\bibitem{grimmett}
G. Grimmett and D. Stirzaker, 
\textit{Probability and Random Processes}, 
4th ed., Oxford University Press, 2020.

\bibitem{schoenfeld}
A. H. Schoenfeld, 
\textit{Mathematical Problem Solving}, 
Academic Press, 1985.

\bibitem{mazur}
B. Mazur, 
\textit{When is One Thing Equal to Some Other Thing?}, 
\textit{The Mathematical Intelligencer}, 12(4), 1990, pp. 58–65.

\bibitem{thom}
R. Thom,
\textit{Modern Mathematics: Does it Exist?},
in: A. G. Howson (ed.), \textit{Developments in Mathematical Education},
Cambridge University Press, 1973, pp.~194--210.

\bibitem{munkres}
J. R. Munkres, 
\textit{Topology}, 2nd ed., 
Prentice Hall, 2000.

\end{thebibliography}
\end{document}